\documentclass[12pt]{amsart}
\usepackage{amssymb}
\usepackage{amscd}

\newtheorem{exercise}{Exercise}[section]

\newcommand{\BN}{{\mathbb{N}}}
\newcommand{\BR}{{\mathbb{R}}}

\newcommand{\BE}{{\mathbb{E}}}
\newcommand{\BF}{{\mathbb{F}}}
\newcommand{\BQ}{{\mathbb{Q}}}
\newcommand{\BP}{{\mathbb{P}}}

\newcommand{\BG}{{\mathbb{G}}}

\newcommand{\BH}{{\mathbb{H}}}
\newcommand{\BK}{{\mathbb{K}}}

\newcommand{\OO}{{\mathcal{O}}}
\newcommand{\gD}{\Delta}
\newcommand{\gL}{\Lambda}
\newcommand{\gd}{\delta}

\newcommand{\gC}{\Gamma}
\newcommand{\gc}{\gamma}
\newcommand{\gs}{\sigma}
\newcommand{\gS}{\Sigma}
\newcommand{\gO}{\Omega}
\newcommand{\gep}{\epsilon}
\newcommand{\gl}{\lambda}
\newcommand{\ga}{\alpha}
\newcommand{\gt}{\tau}

\newcommand{\ti}[1]{\tilde{#1}}

\newcommand{\SL}{\text{SL}}
\newcommand{\GL}{\text{GL}}
\newcommand{\PSL}{\text{PSL}}
\newcommand{\PGL}{\text{PGL}}

\newtheorem{prop}{Proposition}[section]
\newtheorem{thm}[prop]{Theorem}
\newtheorem{lem}[prop]{Lemma}

\newtheorem{cor}[prop]{Corollary}

\theoremstyle{definition}

\newtheorem{defn}[prop]{Definition}
\newtheorem{rem}[prop]{Remark}

\newtheorem{clm}[prop]{Claim}

\begin{document}
\author{E. Breuillard and T. Gelander}

\thanks{E.B. acknowledges support from the French CNRS and the IAS Princeton}
\thanks{T.G. was partially supported by NSF grant
DMS-0404557, and by BSF grant 2004010}

\date{\today}
\title{Uniform independence in linear groups}

\begin{abstract}
We show that for any finitely generated group of matrices that is
not virtually solvable, there is an integer $m$ such that, given
an arbitrary finite generating set for the group, one may find two
elements $a$ and $b$ that are both products of at most $m$
generators, such that $a$ and $b$ are free generators of a free
subgroup. This uniformity result improves the original statement
of the Tits alternative.
\end{abstract}

\maketitle

\setcounter{tocdepth}{1}

\tableofcontents

%------------------------------------------------------------------------------

\section{Introduction}

The main results of this paper were announced in \cite{note}. We
will say that two elements $x,y$ in a group $\gC$ are {\it
independent} if they satisfy no relation, i.e. if they generate a
non-abelian free subgroup. The classical Tits' alternative
\cite{Tits} says that if $\gC$ is a finitely generated linear
group which is not virtually solvable (i.e. does not contain a
solvable subgroup of finite index), then $\gC$ contains two
independent elements. However, Tits' proof gives no indication of
how deep inside the group one has to look in order to find
independent elements. The main result of this paper is the
following:

\begin{thm}\label{thm1}
Let $\gC$ be a finitely generated non-virtually solvable linear
group. Then there is a constant $m=m(\gC )\in\BN$ such that for
any symmetric generating set $\gS$ ($\gS \ni id$) of $\gC$, there
are two words $W_1,W_2$ of length at most $m$ in the alphabet
$\gS$ for which the corresponding elements in $\gC$ are
independent. In other words, the set $\gS^m$ contains two
independent elements.
\end{thm}

By linear group, we mean any subgroup of $\GL_d(K)$ for some
integer $d\geq 1$ and some field $K$. Let now $\BK$ be a global
field, $\overline{\BK}$ its algebraic closure and $S$ a finite set of places of $\BK$ containing all infinite ones.
%, and $\BG$ a simple algebraic group defined over
%$\BK$.
We denote by $\OO_\BK(S)$ the ring of $S$--integers.
% and by
%$\BG (\OO_\BK(S))$ the corresponding $S$--arithmetic group. It is
%well known that $\BG$ has a $\BK$--rational irreducible
%representation $\rho$ into $\SL_d$ with $d=\dim (G)$ such that
%$\rho (\BG (\OO_\BK(S)))\leq\SL_d(\OO_\BK(S))$.
%Since any
%non-virtually solvable linear group has a homomorphism into some
%$S$--arithmetic group with Zariski dense image (see Lemma
%\ref{lem:spec}),
A subgroup of $GL_d(\overline{\BK})$ will be called
\textit{irreducible} if it does not leave invariant any
non-trivial subspace of $\overline{\BK}^d$ (this is sometimes
called absolutely irreducible). After passing to a suitable
homomorphic image (see Lemma \ref{lem:spec}) of the linear group
under consideration, Theorem \ref{thm1} reduces to the following:

\begin{thm}\label{thm2}
Let $\BK$ be a global field, $S$ a finite set of places of $\BK$
containing all the infinite ones, and $d\geq 2$ an integer. Then
there is a constant $m=m(d,\BK,S)$ with the following property.
Suppose that $\gS\subset\SL_d(\OO_\BK (S))$ is a symmetric subset
containing the identity and generating an irreducible subgroup
whose Zariski closure $\BG$ is semisimple and Zariski connected,
then $\gS^m$ contains two independent elements.
\end{thm}

\begin{rem}
In characteristic zero we can actually find two
independent elements in $\gS^m$ which generate a {\it Zariski
dense} subgroup of $\BG$ (see Theorem \ref{thm:Z-dense} and Remark \ref{rem:Z-dense}).
\end{rem}

As in Tits' original proof we use the classical ping-pong lemma
(Lemma \ref{lem:ping-pong}) for the action of the subgroup
generated by $\gS$ on a projective space over some local field.
Since we are in the arithmetic case, there are only finitely many
candidates for the local field, namely the completions $\BK_v$
with $v \in S$. A substantial part of the proof consists in
finding a ``good" metric on the projective space. If $k$ is a
local field and $H\leq\SL_d(k)$ is a semisimple $k$--subgroup with
corresponding symmetric space (or building) $X$, any point in $X$
determines a metric on $k^d$ hence on the projective space $\BP
(k^d)$. For example, the symmetric space $X=SL_d(\BR)/SO_d(\BR)$
is the space of scalar products on $\BR^d$ with a normalized
volume element. Therefore finding a ``good" metric on $\BP (k^d)$
amounts to finding a ``good" point in $X$.

In Section \ref{sec:displacement}, Lemma \ref{prop:constant-c}, we
establish a useful inequality, a norm-versus-spectrum Comparison
Lemma, that relates the displacement of any finite (and more
generally compact) set $\gS$ of isometries of the symmetric space
(or building) of $\SL_d(k)$ to the displacement of a single
element lying in $\gS^{d^2}$. This Comparison Lemma supplies us
with a good metric on $\BP (k^d)$ and an element in $\gS^{d^2}$
that has a ``large" eigenvalue compared to the Lipschitz constants
(for this good metric) of every generator in $\gS$. With such
information, it is not difficult to produce two proximal elements
with distinct attracting points that will generate a free
semi-group. Hence a consequence of our Comparison Lemma is the
Eskin--Mozes--Oh theorem \cite{EMO} on uniform exponential growth
(for details on this implication and improvements in this
direction, see our subsequent paper \cite{entropy}). However, the
Comparison Lemma alone is not sufficient to prove Theorem
\ref{thm1} and produce the required independent elements. As a
matter of fact, it is usually much harder to generate a free
subgroup than a free semi--group.

In Section \ref{sec:arithmetic-displacement} we prove the
following theorem. Let $\BG$ be a semisimple algebraic
$\BK$--subgroup of $\SL_d$. Let $G=\prod_{v\in S}\BG (\BK_v)$
%, where each $\BK_v$ is the completion of $\BK$ with respect to the absolute value $v$.
and $\gC=\BG(\OO_\BK (S))$ be a corresponding $S$--arithmetic
group, which we view as a discrete subgroup of $G$ via the
diagonal embedding. By the Borel Harish-Chandra theorem $\gC$ is a
lattice in $G$, i.e. the quotient space $G/\gC$ carries a finite
$G$--invariant measure. Let $X$ be the product of symmetric spaces
and affine buildings associated to $G$ with a base point $x_0$.

\begin{thm}\label{displacement-growth}
There are positive constants $c_1$ and $c_2$ such that for any
finite subset $\Sigma$ in $\Gamma$ generating a subgroup whose
Zariski closure is connected semisimple and not contained in a
proper parabolic subgroup of $\BG$, we have for all $x \in X$:
\begin{equation*}
d_\gS(x) \geq c_1 d_{X/\gC}(\pi(x),\pi(x_0)) - c_2,
\end{equation*}
where $d_\gS(x)=\max \{ d(g\cdot x,x), g \in \gS \}$, $\pi(x)$ is
the projection of $x$ to the locally symmetric space $X/\gC$ and
$d_{X/\gC}$ is the induced metric on $X/\gC$.
\end{thm}

In other words, the displacement in $X$ of a finite set of lattice
points must grow at a fixed linear rate independently of the
finite set, as one tends to infinity in the locally symmetric
space $X/\gC$, provided that it generates a ``large enough"
subgroup. Note that this theorem is trivial when $\gC$ is uniform.
Moreover, the analogous result holds also for non-arithmetic
lattices (see Remark \ref{rem:n-a-v}) and the constant $c_1$ can
actually be taken to be independent of the choice of the lattice
inside a given group $G$.

At the beginning of the argument proving Theorem
\ref{displacement-growth}, we establish Lemma \ref{k1,l1}, a quantitative version of the Kazhdan--Margulis theorem (namely, if $g\in G$ is "far" from $\gC$ then $g\gC g^{-1}$ contains a non-trivial unipotent ``close" to the identity), which is itself of independent interest.

As a corollary of Theorem \ref{displacement-growth} we obtain
Proposition \ref{prop:gamma}, an arithmetic variant of the
Comparison Lemma. Hence, the outcome of Section
\ref{sec:arithmetic-displacement} is that we can choose the
``good" metric on $\BP (k^d)$ to be arithmetically defined. This
will turn out to be crucial when constructing the ping--pong
players.

Section \ref{sec:ping-pong} is devoted to the construction of the
desired independent elements as ping--pong players on $\BP (k^d)$.
This is done in four steps. First, we construct a proximal
element, second, a very contracting one, third, a very proximal
one, and fourth, a conjugate of the very proximal element which
will form the second ping--pong partner (see Section \ref{prelim}
for this terminology). This construction relies on the study of
the dynamics of projective transformations carried out in
\cite{BG}, and in particular the relation (first used by Tits in
his original proof) between the contraction properties of a
transformation and the Lipschitz constant of its restriction to an
open subset (see Proposition \ref{contracting-properties}). The
arithmetically defined metric that we get from Section
\ref{sec:arithmetic-displacement} supplies us with the two main
ingredients needed to construct the desired ping--pong pair,
namely control on proximality and control on the ability to
separate projective points from projective hyperplanes. The
guiding idea is that the distance between two arithmetically
defined objects is either zero or can be bounded from below by
arithmetic data.

In Section \ref{sec:Z-dense} we restrict to the characteristic $0$ case and show that the bounded independent elements can be chosen to generate a Zariski dense subgroup.

In Section \ref{sec:applications} we describe some consequences of
Theorem \ref{thm1}. One of the main application is that a finitely
generated non--amenable linear group is uniformly non--amenable
and has uniform Cheeger constant, i.e. the family of all Cayley
graphs associated with finite generating sets forms a uniform
family of expanders, see Section \ref{subsec:u-n-a}. One important
consequence is the following:

\begin{thm}\label{thm12345}
Let $\gC$ be a non--virtually solvable linear group. Then there is
a positive constant $\gep$ such that for any finite (not
necessarily symmetric) generating set $\gS$ of $\gC$, and any
finite set $A\subset\gC$, there is some $\gs\in\gS$ for which
$$
 \frac{|\gs A\triangle A|}{|A|}>\gep.
$$
\end{thm}

Theorem \ref{thm12345} has several consequences, for instance, for
the growth function of $\gC$ with respect to a varying generating
set. Clearly it implies that $\gC$ has uniform exponential growth,
but in addition it shows that the growth function gets larger when
the generating set get larger. Moreover since in Theorem
\ref{thm12345} we do not assume, in contrast to the situation in
\cite{EMO}, that the generating set $\gS$ is symmetric, we obtain
a uniform exponential growth result for semi--groups rather than
for groups. As another example, note that it implies uniform
exponential growth for spheres rather than for balls. For more
results in this vein see Section \ref{sub-sec:growth}.

We will also show that Theorem \ref{thm1} implies the connected
case of the Topological Tits Alternative from \cite{BG} and
\cite{BG1}. Recall that the connected case of the Topological Tits
Alternative had several interesting consequences such as the
Connes--Sullivan conjecture about amenable actions of subgroups of
real Lie groups, and the Carri\`{e}re conjecture about the
polynomial versus exponential dichotomy for the growth of leaves
in a Riemannian foliation on a compact manifold. In particular,
Theorem \ref{thm1} implies these results as well, see Section
\ref{seb-sec:dense}.

\medskip{Acknowledgements:}
We thank G.A. Margulis for his interest in this work and for many
conversations and discussions, and in particular for suggesting us
Lemma \ref{k1,l1}. We thank A. Salehi--Golsefidy for many
conversations and suggestions which helped to overcome several
difficulties that arose in the positive characteristic case.

%-------------------------------------------------------------------------------------

\section{Some preliminaries}\label{prelim}

\subsection{Dynamics of projective transformations}
For a more exhaustive and detailed study of the dynamical
properties of projective transformations we refer the reader to
[\cite{BG}, Section 3] and [\cite{BG1}, Section 3].

Let $k$ be a local field and $\left\| \cdot \right\| $ the
standard norm on $ k^{n}$, i.e. the standard Euclidean (resp.
Hermitian) norm when $k$ is ${ \mathbb{R}}$ or ${\mathbb{C}}$ and
$\left\| x\right\| =\max_{1\leq i\leq n}|x_{i}|$ where $x=\sum
x_{i}e_{i}$ when $k$ is non-Archimedean and $ (e_{1},\ldots
,e_{n})$ is the canonical basis of $k^{n}$. This induces an
operator norm on $\text{SL}_{n}(k)$. Consider the standard Cartan
decomposition of $\text{SL}_{n}(k)$,
\begin{equation*}
\text{SL}_{n}(k)=KAK
\end{equation*}
where $K$ is $\text{SO}_{n}(\Bbb{R}),\text{SU}_{n}(\Bbb{C})$ or
$\text{SL}_{n}({\mathcal{O}}_{k})$ according to whether
$k={\mathbb{R}},{\mathbb{C}}$ or is non-Archimedean, and
$A=\{\text{diag}(a_{1},\ldots ,a_{n}):a_{1}\geq \ldots \geq
a_{n}>0,\prod a_{i}=1\}$ if $k$ is Archimedean, and
$A=\{\text{diag} (\pi ^{j_{1}},\ldots ,\pi ^{j_{n}}):j_{i}\in
{\mathbb{Z}},j_{i}\leq j_{i+1},\sum j_{i}=0\}$ if $k$ is
non-Archimedean with uniformizer $\pi $. Any element $g\in
\text{SL}_{n}(k)$ can be decomposed as a product
$g=k_{g}a_{g}k_{g}^{\prime }$, where $k_{g},k_{g}^{\prime }\in K$
and $a_{g}\in A$. The $A$--part $a_{g}$ is uniquely determined by
$g$, but $k_{g},k_{g}^{\prime }$ are not. We will set
\begin{equation*}
a_{g}=\text{diag}(a_{1}(g),\ldots ,a_{n}(g)).
\end{equation*}
Note that $a_{1}(g)=\Vert a(g)\Vert =\Vert g\Vert $. For $g\in
\text{SL} _{n}(k)$ we denote by $[g]$ the corresponding projective
transformation $ [g]\in \text{PSL}_{n}(k)$. Similarly, for $v\in
k^n$ we denote by $[v]$ the corresponding projective point, and
for a linear subspace $H\leq k^n$ we let $[H]$ be the
corresponding projective subspace.

The canonical norm on $k^{n}$ induces the associated canonical
norm on $\bigwedge^{2}k^{n}$. We define the \textit{standard
metric} on $ \Bbb{P}^{n-1}(k)$ by the formula
\begin{equation*}
d\big( [v],[w]\big) =\frac{\left\| v\land w\right\| } {\left\|
v\right\|\cdot\left\| w\right\| }
\end{equation*}
This is well defined and satisfies the following properties:

\textbf{(i)} $d$ is a distance on $\Bbb{P}^{n-1}(k)$ which induces the
canonical topology inherited from the local field $k$.

\textbf{(ii)} $d$ is an ultra--metric distance if $k$ is
non-Archimedean, i.e.
\begin{equation*}
d\big( [v],[w]\big)\leq \max \{d\big( [v],[u]\big) ,d\big( [u],[w]\big)\}
\end{equation*}
for any non-zero vectors $u,v$ and $w$ in $k^{n}$.

\textbf{(iii)} If $f$ is a linear form $k^{n}\to k$, then for any
non-zero vector $v\in k^{n}$,
\begin{equation}
d\big( [v],[\ker f]\big) =\frac{\left| f(v)\right| }{\left\| f\right\| \cdot
\left\| v\right\| }  \label{lin}
\end{equation}

\textbf{(iv)} Every projective transformation
$[g]\in\text{PSL}_n(k)$ is bi--Lipschitz on the entire projective
space with Lipschitz constant $|\frac{a_1(g)}{a_n(g)}|^2=\|
g\|^2\cdot\|g^{-1}\|^2$.

\begin{defn}
A projective transformation $[g]\in \text{PGL}_{n}(k)$ is called
\textit{$\epsilon $--contracting}, for some $\epsilon
>0$, if there is a projective hyperplane $[H]$, called a
\textit{repelling hyperplane}, and a projective point $[v]$,
called an \textit{attracting point} such that for all points $[p]\in \Bbb{P}^{n-1}(k)$,
\begin{equation*}
d([p],[H])\geq \epsilon \Rightarrow d([gp],[v])\leq \epsilon .
\end{equation*}
An element $[g]$ is called \textit{$(r,\epsilon )$--proximal}, for
$ r>2\epsilon $, if it is $\epsilon$--contracting with respect to
some $[H],[v] $ with $d([H],[v])\geq r$. An element $[g]$ is
called \textit{$\epsilon$--very contracting} (resp.
\textit{$(r,\epsilon)$--very proximal}) if both $ [g]$ and
$[g^{-1}]$ are $\epsilon$--contracting (resp.
$(r,\epsilon)$--proximal).
\end{defn}

The following proposition summarizes the relations between
contraction, Lipschitz constants and the ratio between the highest
coefficients of $a_g$.

\begin{prop}[See Lemma 3.4 and 3.5 in \cite{BG} and Proposition 3.3 in \cite{BG1}]
\label{contracting-properties} Let $\epsilon \in
(0,\frac{1}{4}],~r\in (0,1]$. Let $g\in \text{SL}_{n}(k)$.

\begin{enumerate}
\item If $|a_{2}(g)/a_{1}(g)|\leq \epsilon $ then $[g]$ is
$\epsilon /r^{2}$--Lipschitz outside the $r$--neighborhood of the
repelling hyperplane $[ \text{span}\{k^{\prime
}{}_{g}^{-1}(e_{i})\}_{i=1}^{n}]$.

\item If the restriction of $[g]$ to some open subset $O\subset {
\mathbb{P}}^{n-1}(k)$ is $\epsilon $--Lipschitz, then $
|a_{2}(g)/a_{1}(g)|\leq \epsilon /2$.

\item If $|a_{2}(g)/a_{1}(g)|\leq \epsilon ^{2}$ then $[g]$ is
$\epsilon $ -contracting, and vice versa, if $[g]$ is $\epsilon
$--contracting, then $ |a_{2}(g)/a_{1}(g)|\leq c\epsilon ^{2}$
where $c$ is some constant depending on $k$.
\end{enumerate}
\end{prop}

Note that the attracting point and repelling hyperplane of a
contracting or proximal element are not uniquely defined. In case
$g$ is semisimple, it is sometimes useful to choose them to be the
span of relevant eigenvectors of $g$, while it is also possible to
define them using the Cartan decomposition like in point $(1)$
above. Very proximal elements are our tool to generate free
subgroups via the following version of the classical ping-pong
lemma:

\begin{lem}[The Ping--Pong Lemma]\label{lem:ping-pong}
Assume that $x$ and $y$ are $(r,\epsilon )$--very proximal
projective transformations of $ \Bbb{P}^{n-1}(k)$ (for some
$r>2\epsilon $), and suppose that the distances between the
attracting points of $x^{\pm 1}$ (resp. of $y^{\pm 1}$) and the
repelling hyperplanes of $y^{\pm 1}$ (resp. of $x^{\pm 1}$) are at
least $r$, then $x$ and $y$ are independent.
\end{lem}

%----------------------------------------------------------------

\subsection{How to get out of proper subvarieties in bounded time}
\label{subsection:verities}
\addcontentsline{toc}{subsection}{Bezout's theorem}

Recall the following classical theorem (c.f. \cite{Sch}):

\begin{thm}[Generalized Bezout theorem]\label{thm:Bezout}
Let $\BK$ be a field, and let $X_1,\ldots,X_s$ be pure dimensional
algebraic subvarieties of $\BK^n$. Denote by $Z_1,\ldots,Z_t$ the
irreducible components of $X_1\cap\ldots\cap X_s$. Then
$$
 \sum_{i=1}^t \text{deg}(Z_i)\leq \prod_{j=1}^s\text{deg}(X_j).
$$
\end{thm}

For an algebraic variety $X$ we will denote by $\chi (X)$ the sum
of the degrees and dimensions of its irreducible components. The following lemma
is a consequence of Theorem \ref{thm:Bezout} (see Lemma 3.2 in
\cite{EMO} and its proof\footnote{In \cite{EMO} it is assumed that
$\Sigma$ is finite, that the characteristic of the field is $0$,
and that the algebraic group $\BG$ and the variety $X$ are fixed, however the proof in \cite{EMO}
does not depend on these assumptions.}).

\begin{lem}\cite{EMO}\label{Bezout}
Given an integer $\chi$ there is $N=N(\chi)$ such that for any
field $K$, any integer $d\geq 1$, any $K$--algebraic subvariety
$X$ in $GL_d(K)$ with $\chi (X)\leq\chi$ and any subset $\Sigma
\subset {GL_d(K)}$ which contains the identity and generates a
subgroup which is not contained in $X(K)$, we have $\Sigma ^{N}
\nsubseteq X(K)$.
\end{lem}

When $X$ is given, we will sometimes abuse notations and
write $N(X)$ for $N(\chi (X))$.

%-----------------------------------------------------

\section{Reduction to the $S$--arithmetic setting}

Here we reduce Theorem \ref{thm1} to Theorem \ref{thm2}. Given a
global field $\BK$ and a finite set $S$ of places of $\BK$
including all the infinite ones, we denote by $\OO_\BK (S)$ the
ring of $S$--integers in $\BK$. The following lemma is well known:

\begin{lem}\label{lem:spec}
Let $\gC$ be a finitely generated linear group which is not
virtually solvable. Then there is a global field $\BK$, a finite
set of places $S$ of $\BK$ and a representation
$f:\gC'\to\GL_d(\OO_\BK (S))$ of some finite index subgroup
$\gC'\leq\gC$ whose image is Zariski dense in a simple
$\BK$--algebraic group.
\end{lem}

\begin{proof}(Suggested to us by G.A. Margulis)
In the proof of the classical Tits alternative \cite{Tits}, Tits
produces a local field $k$ and a homomorphism $\varphi
:\Gamma\rightarrow\text{GL}_n(k)$ such that $\varphi(\Gamma)$
contains two proximal elements $\varphi
(x),\varphi(y)$ which are ``playing ping--pong" on the projective space ${\mathbb{P}}%
^{n-1}(k)$ (i.e. satisfy the hypothesis of Lemma \ref{lem:ping-pong}) and hence generate a free subgroup.

Let $F$ be a global field whose completion is $k$, and let
$\overline{F}$ be its integral closure in $k$, i.e. the field of
all elements in $k$ algebraic over $F$. Let
$X=\text{Hom}(\Gamma,\text{GL}_n(k))$ be the variety of
representations of $\Gamma$ into $\text{GL}_n(k)$, realized as a subset of $%
\text{GL}_n(k)^{d(\Gamma )}$ where $d(\gC )$ is the size of some
finite generating set of $\gC$. Then $X$ is an algebraic variety
defined over $F$, and as follows from the implicit function
theorem, the set $X(\overline{F })$ of $\overline{F}$ points is
dense in $X(k)$ in the topology
induced from $\text{GL}_n(k)^{d(\Gamma )}$. Thus we can choose a deformation $%
\rho\in X(\overline{F})$ arbitrarily close to $\varphi$. Now if
$\rho$ is sufficiently close to $\varphi$, then $\rho (x)$ and
$\rho (y)$ still play ping--pong on ${\mathbb{P}}^{n-1}(k)$, and
this implies that $\rho (\Gamma )$ is not virtually solvable.

Let ${\mathbb{K}}$ be the field generated by the entries of
$\rho (\Gamma )$. Since $\Gamma$ is finitely generated,
${\mathbb{K}}$ is a global field. Let $\Gamma'$ be a finite
index subgroup of $\gC$ such that $\rho (\Gamma')$ is Zariski connected. We then obtain the representation $%
f$ and the group ${\mathbb{G}}$ by dividing by the solvable
radical and projecting to a simple factor of the Zariski closure
of $\rho(\Gamma')$. Note that as $\rho
(\Gamma')\subset\text{GL}_n({\mathbb{K}} )$ its Zariski
closure and solvable radical are defined over ${\mathbb{K}}$. Therefore $%
f(\Gamma')\leq{\mathbb{G}} ({\mathbb{K}} )$.

Finally, since $\Gamma$ is finitely generated, there is a finite
set of
places $S$ such that $f(\Gamma )$ lies in the $S$-arithmetic group ${%
\mathbb{G}} ({\mathcal{O}}_\BK (S))$.
\end{proof}

It is easy to check that if $n$ is the index of $\gC'$ inside
$\gC$, then for any generating set $\gS\ni 1$ of $\gC$ containing
the identity, $\gS^{2n+1}$ contains a generating set for $\gC'$.
Hence Lemma \ref{lem:spec} implies that Theorem \ref{thm1} is a
consequence of Theorem \ref{thm2}. The main part of this paper is
therefore devoted to the proof of \ref{thm2}.

%------------------------------------------------------------------------------

\section{Minimal norm versus Maximal eigenvalue}
\label{sec:displacement}

In this section, we state and prove the Comparison Lemma, Lemma
\ref{prop:constant-c}. Roughly speaking, this statement says that
the minimal displacement of a compact subset of isometries of a
symmetric space or an affine building is comparable to the minimal
displacement of a single element belonging to some bounded power
of the subset. When we came up with Lemma \ref {prop:constant-c},
we were strongly inspired by Proposition 8.5 in \cite{EMO}.

%-------------

\subsection{Minimal norm, maximal eigenvalue, and the Comparison Lemma}

Let $k$ be a local field with absolute value $|\cdot |_{k}.$ It
induces the standard norm on $k^{d}$ which in turn gives rise to
an operator norm $\left\| \cdot \right\| $ on $M_{d}(k).$ If $k$
is not Archimedean, let $\mathcal{O}_{k}$ be its ring of integers
and $m_{k}$ the maximal ideal in $\mathcal{O}_{k}.$ We note that
$\left\| a\right\| _{k}\geq 1$ for all $a\in\SL_{d}(k).$ Let
$\Lambda _{k}(a)$ be the maximum absolute value of all eigenvalues
of $a$ (recall that the absolute value has a unique extension to
the algebraic closure of $k$). If $g\in\SL_{d}(k)$ we denote by
$a^{g}$ the conjugate $gag^{-1}$.

For a compact subset $Q\subset M_{d}(k)$ we denote:
\begin{eqnarray*}
\Lambda _{k}(Q) &=&\max \{\Lambda _{k}(a):a\in Q\} \\
\Vert Q\Vert _{k} &=&\max \{\Vert a\Vert _{k}:a\in Q\}~ \\
\gD_{k}(Q) &=&\inf_{g\in \text{SL}_{d}(k)}\Vert gQg^{-1}\Vert
_{k}.
\end{eqnarray*}

\begin{rem}
One can define $\ti \gD_k$ by taking the infimum over
$g\in\PGL_d(k)$. This has some advantages in the non-Archimedean
case, e.g. $\ti\gD_k(Q)=1$ whenever $Q$ lies in a compact group.
Moreover, the ratio between $\ti\gD_k$ and $\gD_k$ is bounded
since $\PSL_d$ has finite index in $\PGL_d$. However, we found it
more convenient for our purposes to use $\gD_k$ as defined above.
\end{rem}

%Note that $\Lambda _{k}(Q)^{n}\leq \Lambda _{k}(Q^{n})\leq
%\gD_{k}(Q^{n})\leq \gD_{k}(Q)^{n}$ for all $n\geq 1$, where $Q^n$
%is the set of all products of $n$, not necessarily different,
%elements of $Q$.
In terms of the action of $\text{SL}_{d}(k)$ on its symmetric
space or affine building, $\log \gD_{k}(Q)$ is, up to a
multiplicative constant, the minimal displacement of $Q$, i.e. the
smallest radius of a $Q$--orbit. When $Q=\{a\}$ is a single
element, diagonalizable over $k$, we have $\gD_{k}(\{a\})=\Lambda
_{k}(a)$. The following gives a similar relation when $Q$ is an
arbitrary compact subset. We denote by $Q^i$ the set of all
products of $i$, not necessarily different, elements of $Q$.

\begin{lem} \label{prop:constant-c} \textbf{(Norm--versus--Spectrum Comparison Lemma)}
There exists a constant $c=c(d,k)>0$ such that for any
compact subset $Q\subset M_{d}(k)$ we have
\begin{equation*}
\gD_{k}(Q)^{i}\geq \Lambda _{k}(Q^{i})\geq c\cdot \gD_{k}(Q)^{i}
\end{equation*}
for some $i\leq d^{2}$.
%Moreover, if $k$ is non-Archimedean the
%inequality holds with $c=\left( |\pi |_{k}\right) ^{2d-1}$ for a
%uniformizer $\pi $ for $k$.
\end{lem}

\begin{rem}
The proof that we are about to give uses a compactness argument
and hence is not effective. In \cite{entropy} we will give an
effective proof of \ref{prop:constant-c}. This relies on an
effective proof of Wedderburn's theorem on the existence of
idempotents in non-nilpotent subalgebras of matrices.
Additionally, we will show in \cite{entropy}  that when $k$ is
non-Archimedean, by taking finite extensions, we can make $c$
arbitrarily close to $1$ (actually $c=(|\pi |_{k})^{2d-1}$), and
derive a strong uniformity result concerning the growth functions
of linear groups.
\end{rem}

We now proceed to the proof of Lemma \ref{prop:constant-c}. We
start with the following classical statement:

\begin{lem}
\label{nilp:span} Let $R$ be a field or a finite ring and let $\mathcal{A}%
\leq M_{d}(R)$ be a subring and $R$--submodule. Suppose that
$\mathcal{A}$ is spanned as an $R$--module by nilpotent matrices,
then $\mathcal{A}$ is nilpotent, i.e. $\mathcal{A}^{n}=\{0\}$ for
some $n\geq 1.$
\end{lem}

\begin{proof}
In the $0$ characteristic case, the lemma follows easily from
Engel's theorem using the fact that a matrix is nilpotent iff all
its powers have $0$ trace. The proof we give now works in arbitrary characteristic and
 was suggested to us by A. Salehi-Golsefidy.

The ring $\mathcal{A}$ is Artinian and
therefore its Jacobson radical
%\footnote{Recall that
%$J(\mathcal{A})$ is the intersection of the annihilators of all simple $\mathcal{A}-\text{modules}\}$.}
$J(\mathcal{A})$ is nilpotent. We will prove the lemma by showing
that $\mathcal{A}=J(\mathcal{A})$. Let
$\mathcal{B}=\mathcal{A}/J(\mathcal{A})$ and assume by way of
contradiction that $\mathcal{B}\neq 0$. Now $\mathcal{B}$ is
semisimple, hence by the Artin--Wedderburn theorem,
$\mathcal{B}\cong\bigoplus M_{d_i}(\mathcal{D}_i)$, where the
$\mathcal{D}_i$ are division rings. Since $\mathcal{A}$ is spanned
by nilpotent elements, so is $\mathcal{B}$. This implies that the
trace of any element in $M_{d_i}(\overline{k})$ is $0$, and hence
that $d_i=0$. A contradiction.
\end{proof}

Note that an element $A\in M_{d}(k)$ is nilpotent iff $\Lambda_{k}(A)=0$.
The following generalizes this statement to compact subsets.

\begin{lem}
\label{Lambda&Delta=0} For a compact subset $Q\subset M_{d}(k)$ the
following are equivalent:

$(i)$ $Q$ generates a nilpotent subalgebra.

$(ii)$ $\gD_{k}(Q)=0$.

$(iii)$ $\Lambda _{k}(Q^{i})=0,~\forall i\leq d^{2}$.
\end{lem}

\begin{proof}
Let $\mathcal{A}$ be the algebra generated by $Q$.

\noindent $(i)\Rightarrow (ii)$: By Engel's theorem $\mathcal{A}$
and hence $Q$ can be conjugated by a matrix in SL$_{d}(k)$ into
the algebra of upper triangular matrices with $0$ diagonal.
Conjugating further by some suitable diagonal element in
SL$_{d}(k)$ we can make the norm of $Q$ arbitrarily small.

\noindent $(ii)\Rightarrow (iii)$: For any element $g\in M_d
(k),~\|g\|\geq\Lambda_k(g)$, hence $\gD_k(Q)^i\geq
\gD_k(Q^i)\geq\Lambda_k(Q^i)$.

\noindent $(iii)\Rightarrow (i)$: Take $q\leq d^2$ such that $\dim
(\text{span}\cup_{j=1}^q Q^j)=\dim (\text{span}\cup_{j=1}^{q+1}
Q^{j})$, then $\cup_{j=1}^q Q^j$ spans $\mathcal{A}$. Since
$\Lambda (Q^i)=0$ for $i\leq d^2$, it consists of nilpotent
elements; hence the implication follows from Lemma
\ref{nilp:span}.
\end{proof}

\begin{proof}[Proof of Lemma \ref{prop:constant-c}]
Suppose by contradiction that there is a sequence of compact sets
$Q_1,Q_2,\ldots$ in $M_d(k)$ such that $\Lambda_k(Q_n^i)<
{\gD_k(Q_n^i)}/{n},~\forall i\leq d^2$. By replacing $Q_n$ with a
suitable conjugate of it, we may assume that $\| Q_n\|_k\leq
2\gD_k (Q_n)$, and by normalizing we may assume that $\|
Q_n\|_k=1$. Let $Q$ be a limit of $Q_n$ with respect to the
Hausdorff topology on $M_d(k)$. Then $\| Q\|_k=1$, $\gD_k (Q)\geq
\frac{1}{2}$ since $\gD_k$ is upper semi-continuous, and by
continuity of $\Lambda_k$, $\Lambda_k(Q^i)=0,~\forall i\leq d^2$.
This however contradicts Lemma \ref{Lambda&Delta=0}.
\end{proof}

%---------------------

\subsection{Geometric interpretation of the Comparison Lemma}

For $g\in \text{SL}_{d}(k)$ and $x$ in the associated symmetric
space (resp. affine building) $X$, we denote by $d_{g}(x)=d(g\cdot
x,x)$ the displacement of $g$
at $x$. Similarly, for a compact set $Q\subset \text{SL}_{d}(k)$ we let $%
d_{Q}(x)=\max_{g\in Q}d_{g}(x)$. Finally, we consider the minimal
displacement of $g,$ or $Q,$ namely $d_{g}:=\inf_{x\in X}d_{g}(x)$ and $%
d_{Q}:=\inf_{x\in X}d_{Q}(x).$

Therefore, Lemma \ref{prop:constant-c} implies the following
geometric statement:

\begin{cor}\label{GeomInt}
There is a universal constant $C=C(d)>0$ such that for any compact set $%
Q\subset \text{SL}_{d}(k)$ there exists $g\in \cup _{1\leq i\leq d^{2}}Q^{i}$
such that
\begin{equation*}
\frac{1}{\sqrt{d}}d_{Q}-C\leq d_{g}\leq d^{2}\cdot d_{Q}
\end{equation*}
\end{cor}

\proof Clearly, if $g\in Q^{i}$, $d_{g}\leq d_{Q^{i}}\leq i\cdot
d_{Q}.$ Note that (see Lemma \ref{lem:comp} below) for every $x\in
X,$ and $g\in $SL$_{d}(k),$ we have $\log \left\| g\right\|
_{x}\leq d_{g}(x)\leq \sqrt{d}\log \left\| g\right\| _{x}$ where
$\left\| \cdot \right\| _{x}$ is the norm associated to the
compact stabilizer of $x$ in SL$_{d}(k)$, and the $\log $ is taken
in base $|\pi |_{k}^{-1}$ when $k$ is non-Archimedean. Since the
action of $\SL_d(k)$ on $X$ is transitive in the Archimedean case
and transitive on the cells in the non-Archimedean case, it
follows that $\log \gD_{k}(Q)\geq \frac{1}{\sqrt{d}}(d_{Q}-2)$. On
the other hand, $\log \Lambda _{k}(g)\leq d_{g}$ for all $g\in
$SL$_{d}(k),$ and by Lemma \ref
{prop:constant-c}, there exists an $i\leq d^{2}$ and $g\in Q^{i}$ with $%
\Lambda _{k}(g)\geq c\cdot \gD_{k}(Q)^{i}.$ Hence $d_{g}\geq \log
\Lambda _{k}(g)\geq \frac{i}{\sqrt{d}}(d_{Q}-2)+\log c$.
\endproof

%--------------------------------------------------------------------------------------

\section{Uniform linear growth of displacement functions}\label{sec:arithmetic-displacement}

In this section we prove Theorem \ref{displacement-growth} and
derive an arithmetic analog to Lemma \ref{prop:constant-c} that
will be crucial in the proof of Theorem \ref{thm2}.

Let ${\mathbb{K}}$ be a global field, $S$ a finite set of places
containing all infinite ones and ${\mathcal{O}}_\BK (S)$ the ring of $%
S$--integers. For $v\in S$ we let ${\mathbb{K}}_v$ denote the completion of ${%
\mathbb{K}}$ with respect to $v$. Since $v$ extends uniquely to
any finite extension of ${\mathbb{K}}_v$ we will, abusing notations, denote by $%
|\cdot |_v$ also the absolute value on any such extension. Let
${\mathbb{G}}\leq\SL_d$ be a Zariski connected semisimple
${\mathbb{K}}$--algebraic group.
%Suppose $%
%d=\dim ({\mathbb{G}} )$ and identify ${\mathbb{G}}$ with its image under the
%adjoint representation. Then ${\mathbb{G}}$ is a closed ${\mathbb{K}}$%
%-subgroup of ${\mathbb{SL}}_d$.
Let
\begin{equation*}
G=\prod_{v\in S}{\mathbb{G}} ({\mathbb{K}}_v)\leq H=\prod_{v\in S}{\SL_d ({\mathbb{K}}_v)}.
\end{equation*}
The group of $S$--integers ${\mathbb{G}} ({\mathcal{O}}_\BK (S))$
is an $S$--arithmetic group. We will identify it with its diagonal
embedding in $G$. This makes ${\mathbb{G}} ({\mathcal{O}}_\BK
(S))$ a discrete subgroup of $G$. The Borel Harish-Chandra theorem
says that it is a lattice in $G$, i.e. the quotient space
$G/{\mathbb{G}}({\mathcal{O}}_\BK(S))$ carries a finite
$G$--invariant measure, and that if ${\mathbb{G}}$ is
${\mathbb{K}}$--anisotropic then
$G/{\mathbb{G}}({\mathcal{O}}_\BK(S))$ is compact. We will set
$\Gamma ={\mathbb{G}}({\mathcal{O}}_\BK(S))$.

Consider $v\in S$ and set $G_v={\mathbb{G}}
({\mathbb{K}}_v )$ and $H_v=\text{SL}_d({%
\mathbb{K}}_v)$. Let $K_v\leq\text{SL}_d({%
\mathbb{K}}_v)$ be the maximal compact subgroup corresponding to
the standard norm on ${\mathbb{K}}_v^d$. Recall that for $v$
Archimedean any maximal compact is conjugate to $K_v$ in
$\text{SL}_d({\mathbb{K}}_v)$, and for $v$ non-Archimedean there
are $d+1$ conjugacy classes. Let
$X_v=\text{SL}_d({\mathbb{K}}_v)/K_v$ be the
associated symmetric space or affine building, let $A_v$ be a Cartan semigroup of $\text{SL}_d({%
\mathbb{K}}_v)$ corresponding to $K_v$ with respect to a Cartan
decomposition of $\text{SL}_d({\mathbb{K}}_v)$, and let $x_0\in
X_v$ be the point corresponding to $K_v$.

We also set the following notations. For $a\in \SL_d({\mathbb{K}})$ let
\begin{equation*}
\Lambda (a)=\max \{|\lambda |_{v}:\lambda \text{~is an eigenvalue of~}%
a,~v\in S\}=\max \{\Lambda _{v}(a):v\in S\}.
\end{equation*}
For $v\in S$ let $\Vert \cdot \Vert _{v}$ be the standard operator norm on $\text{SL}_{d}(%
{\mathbb{K}}_{v})$, and for $g=(g_{v})_{v\in S}\in H$ let
\begin{equation*}
\Vert g\Vert =\max \{\Vert g_{v}\Vert _{v}:v\in S\}.
\end{equation*}
For a compact subset $Q\subset H$ we let
\begin{equation*}
\Vert Q\Vert =\max_{a\in Q}\Vert a\Vert ,~~\gD(Q)=\inf_{h\in
H}\Vert Q^{h}\Vert=\max_{v\in S}\Delta_{\BK_{v}}(Q)
,~~\text{and}~~\Lambda (Q)=\max_{a\in Q}\Lambda (a)=\max_{v\in
S}\Lambda_{\BK_{v}}(Q).
\end{equation*}

%------------------------------------------------------------

\subsection{Restating Theorem \ref{displacement-growth}}
\begin{defn}
We will say that a subgroup of $\BG$ is irreducible in $\BG$ if
it is not contained in a proper parabolic subgroup of $\BG$.
\end{defn}

Recall the following result of Mostow in the Archimedean case
(c.f. \cite {Pla-Rap} Theorem 3.7) and Landvogt in the
non-Archimedean (c.f. \cite{Landvogt}):

\begin{thm}
\label{thm:ML} There exists a point $x_1\in X_v$ when $v$ is
Archimedean (resp. a cell $\gs_1\subset X_v$ when $v$ is
non-Archimedean) such that the orbit $G_v\cdot x_0$ (resp.
$\cup\{g\cdot\gs:{g\in G_v}\}$) is convex and isometric to the
symmetric space (resp. affine building) associated to $G_v$.
\end{thm}

Recall that the norm of a matrix in $\text{SL}_d$ is comparable to
the exponent of its displacement. More precisely:

\begin{lem}
\label{lem:comp}For any $h\in \text{SL}_{d}({\mathbb{K}}_{v})$ we
have:

\begin{itemize}
\item  $\Vert h\Vert \leq e^{d(h\cdot x_{0},x_{0})}\leq \Vert h\Vert ^{\sqrt{%
d}}$.

\item  If $x=g^{-1}\cdot x_{0}$ then $\Vert h^{g}\Vert \leq
e^{d(h\cdot x,x)}\leq \Vert h^{g}\Vert ^{\sqrt{d}}$.
\end{itemize}
\end{lem}

\begin{proof}
If $h =k_h a_h k_h'$ is a $KAK$ expression for $h$ then $\| h\|
=\| a_h\|$ and $d(h\cdot x_0,x_0)=d(a_h\cdot x_0,x_0)$ hence its
enough to prove the first inequality for elements in $A$, and for
such elements it follows by a direct computation.

The second inequality is a direct consequence of the first one.
\end{proof}

Assume that $\BG$ is isotropic over $\BK$, i.e that $\Gamma$ is a
non-uniform lattice in $G$. Let $\pi :G\to G/\Gamma$ be the
canonical projection, and for $g\in G$ denote
\begin{equation*}
\|\pi (g)\|=\min_{\gamma\in\Gamma}\| g\gamma\|.
\end{equation*}

Note that the convex orbit of $G$ from Theorem \ref{thm:ML} may
not pass through the origin $x_0$, however, since any two orbits
of $G$ are equidistant, in view of Lemma \ref{lem:comp}, Theorem
\ref{displacement-growth} can be restated as follows:

\begin{thm}\label{thm:DG}
There are positive constants $C_1,C_2$ such that for any finite subset $\Sigma$
in $\Gamma$ generating a subgroup whose Zariski closure is
semisimple and irreducible in $\BG$, we have $\forall g\in G$
\begin{equation}
  \|\gS^g\|\geq C_2\|\pi(g)\|^{C_1}.
\end{equation}
\end{thm}

\begin{rem}\label{rem:n-a-v}
The statement of Theorem \ref{thm:DG}, as well as of Lemma
\ref{k1,l1} below, remains true without the assumption that the
non-uniform lattice $\gC$ is arithmetic.
% but then the constant on
%the exponent $C_G$ has to be replaced by a constant which depends
%on $\gC$, rather than merely on $G$.
To see this one carries the
same argument as below, using a variant of Corollary 8.16 from
\cite{raghunathan} instead of Lemma \ref{epsilon}. Moreover, the constant $C_1$ can be taken to depend only on $G$
and not on the choice of the lattice $\gC$.
\end{rem}

%--------------------------------------------------------

\subsection{A quantitative Kazhdan--Margulis Theorem }
Let $\BK,S,\BG,G,\gC$ be as in the previous paragraph, in
particular we assume that $G/\gC$ is non-compact.

According to the Kazhdan--Margulis Theorem (see
\cite{raghunathan}), if $\|\pi (g)\|$ is large enough, then
$\Gamma ^{g}$ contains a non-trivial unipotent close to the
identity. The following quantitative version of this theorem was
suggested to us by G.A. Margulis.

\begin{lem}
\label{k1,l1} There are positive constants $k_{\gC},l_{\Gamma }$
such for any $g\in
G$ the lattice $\Gamma ^{g}=g\gC g^{-1}$ contains a non-trivial
unipotent $u\in \Gamma ^{g}$ with
\begin{equation*}
\Vert u-1\Vert \leq l_{\Gamma }\Vert \pi (g)\Vert ^{-k_{\gC}}.
\end{equation*}
\end{lem}

Lemma \ref{k1,l1} is proved along the same lines as the original
Kazhdan-Margulis Theorem.

\begin{lem}
\label{epsilon} There is a positive constant $\epsilon _{G}$ such that if $u_{1},\ldots ,u_{t}$ are elements
belonging to a non-uniform $S$--arithmetic subgroup of $G$ and
$\Vert u_{i}-1\Vert \leq \epsilon _{G},~\forall i\leq t$, then the
group $\langle u_{1},\ldots ,u_{t}\rangle $ is unipotent.
\end{lem}

\begin{proof}[Proof of Lemma \ref{epsilon}]
If $\gep_G$ is sufficiently small then by the Zassenhaus Lemma
(c.f. \cite{raghunathan} 8.8. and 8.17.) the $u_i$'s generate a
nilpotent group, and by \cite{margulis} 4.21(A) the $u_i$ are
unipotent. The result follows since any nilpotent group which is
generated by unipotent elements is unipotent.
\end{proof}

\begin{proof}[Proof of Lemma \ref{k1,l1}]
We will first assume that $\text{char}({\BK})=0$, and later
indicate the changes to be made in the positive characteristic
case.

For any Zariski connected unipotent group $U$ there is an element
$g_U\in G$ such that the restriction of $\text{Ad}(g_U)$ to the
Lie algebra of $U$ expands the norm of any element by at least a
factor $4$. Since the Grassmann manifolds are compact, it follows
that there are finitely many elements $g_1,\ldots
,g_k,~g_i=g_{U_i}$ such that for any Lie subalgebra
$\mathfrak{u}$, corresponding to some unipotent subgroup, there is
$i\leq k$ such that the restriction of $\text{Ad}(g_i)$ to
$\mathfrak{u}$ expands the norm of any element by at least a
factor $3$. Now since the exponential map $\exp :\text{Lie}(G)\to
G$ is a diffeomorphism near the origin $0$ of Lie$(G)$ with
differential $1$ at $0$ it follows that for some $\gep_1>0$,
smaller than $\gep_G$, we have:

\begin{equation}
 \|g_iug_i^{-1}-1\|\geq 2\| u-1\|,~~\forall~u\in\exp (\mathfrak{u})
 ~\text{with}~\|u-1\|\leq\gep_1.
\label{star}
\end{equation}
Fix
$$
 a=\max_{i\leq k}\| g_i^{\pm 1}\|.
$$
and let
$$
 k_{\gC}=\log_a 2.
$$
Fix $\gep_2>0$ smaller than $\gep_1$, such that
$$
 B_{\gep_2}(1_G)\subset\cap_{i\leq k}(B_{\gep_1}(1_G))^{g_i^{\pm
 1}}.
$$
Since $M=G/\gC$ has finite volume the ``$\gep_2$--thick part"
$$
 M_{\geq\gep_2}:=\{\pi (g):g\in G,~\text{and}~\gC^g\cap
 B_{\gep_2}(1_G)=\{ 1\}\}
$$
is compact. Let
$$
 l_\gC=\sup_{\{ g: \pi (g)\in M_{\geq\gep_2}\}}\big(\min_{\text{all
 unipotents}~u\in\gC\setminus\{ 1\}}\| u^g-1\|\big)\cdot \sup_{\{ g: \pi (g)\in
 M_{\geq\gep_2}\}}(\|\pi (g)\|)^{k_{\gC}},
$$
then the lemma holds for any $g$ with $\pi (g)\in M_{\geq\gep_2}$.

Now suppose $\pi (g)\notin M_{\geq\gep_2}$. Then $\gC^g$ has a
non-trivial unipotent in $B_{\gep_2}(1_G)$. Let
$$
 b=\min \{\| u^g-1\| :u\in\gC\setminus\{ 1\}~\text{unipotent}\}.
$$
By Lemma \ref{epsilon} $\gC^g\cap B_{\gep_1}(1)$ is contained in
Zariski connected unipotent group, and hence by $(\ref{star})$
there is some $g_{i_1}$ such that the conjugation by it increases
the distance of any non-trivial element of this intersection by at
least a factor of $2$. After this conjugation, there might be some
new unipotent elements in the $\gep_1$--ball around $1_G$,
however, by the choice of $\gep_2$ there are no new unipotents in
the $\gep_2$--ball. Therefore we can iterate this argument
$\lceil\log_2\frac{\gep_2}{b}\rceil$ times, and get a sequence
$g_{i_1},\ldots,g_{i_t},~t={\lceil\log_2\frac{\gep_2}{b}\rceil}$,
such that $\gC^{g_{i_{t}}\cdot\ldots\cdot g_{i_1}g}$ intersect
$B_{\gep_2}(1_G)$ trivially. It follows that $\pi
(g_{i_{t}}\cdot\ldots\cdot g_{i_1}g)\in M_{\geq\gep_2}$, and
hence, if $u\in\gC^{g_{i_{t}}\cdot\ldots\cdot g_{i_1}g}$ is a
non-trivial unipotent closest to $1_G$
$$
 \|\pi (g_{i_{t}}\cdot\ldots\cdot g_{i_1}g)\|^{k_{\gC}}\cdot\|
 u-1\|\leq l_\gC.
$$
Since $\| u-1\|\geq 2^tb$, and since all the $g_i$'s have norm at
most $a$, the result follows.

\medskip

Let us now explain the required modifications in the proof for the
positive characteristic case. For the positive characteristic
version of the Kazhdan--Margulis theorem see \cite{raghunathan1}.
In the positive characteristic case, the unipotent group provided
by Lemma \ref{epsilon} is not Zariski connected, in fact it is
finite. However, it was shown by Borel and Tits \cite{borel-tits}
that for any unipotent group $U$ there is a canonical parabolic
group $P(U)$ which contains the normalizer of $U$ and contains $U$
in its unipotent radical. The unipotent radical of a parabolic
subgroup is Zariski connected, and pro-$p$. Using the $KP$
decomposition where $P$ is a minimal parabolic and $K$ is a
maximal compact, it is easy to show that there is some compact set
$C$ such that for any parabolic subgroup there is $g\in C$ such
that conjugation by $g$ expends the norm of each element in the
unipotent radical of the parabolic by at least $4$, and one can
carry out the same argument as above.
\end{proof}

\subsection{Proof of Theorem \ref{thm:DG}}
Let $\gC,\BG$, $G$, $k_{\gC}$ and $l_{\gC}$ be as in the previous paragraph. Clearly, the
following claim implies Theorem \ref{thm:DG}:

\begin{clm}\label{lem:eq3}
There is a constant $N$, depending only on $\BG$, such that
$\forall g\in G$
\begin{equation}
l_\gC\|\pi(g)\|^{-k_{\gC}}\|\Sigma^g\|^{2N}\geq\epsilon_G.
\label{estimate pi(g)}
\end{equation}
\end{clm}

\begin{proof}
Assume first that char$({\mathbb{K}} )=0$ and let $N=d^2$. Suppose by
way of contradiction that the lemma is false, and let
$u\in\Gamma^g\setminus\{ 1\}$ be a unipotent element as in Lemma
\ref{k1,l1} with
\begin{equation*}\label{eq18}
\|u-1\|\leq l_\gC\|\pi (g)\|^{-k_{\gC}}.
\end{equation*}
Then it follows that for any word $W$ in the elements of
$\Sigma^g$ of length at most $d^2$ we have $\|
u^W-1\|\le\epsilon_G$. Let $\mathcal{U}_i$ be the Zariski closure
of the group generated by $\{ u^W:W$ is a word in the
elements of $\Sigma^g$ of length $\leq i\}$. Then by Lemma \ref{epsilon}, $%
\mathcal{U}_i$ is a unipotent group, hence is Zariski connected.
Therefore,
for some $i_0<\dim ({\mathbb{G}} )\leq d^2$ we have $\mathcal{U}_{i_0}=%
\mathcal{U}_{i_0+1}$, and hence $\mathcal{U}_{i_0}$ is normalized by $%
\Sigma^g$. But this implies that $\Sigma^g$ is contained in some
proper parabolic subgroup, a contradiction to the assumption that
$\Sigma$ generates an irreducible subgroup. Hence the claim is
proved.

%We now give an alternative proof which holds in arbitrary
%characteristic: Let $U$ be a maximal unipotent subgroup of $\BG$,
%and let $X=\{ h\in\BG:U^h\cap U\neq 1\}$. By the Bruhat
%decomposition theorem $X$ is a proper algebraic subvariety of $G$.
%Let
%$$
% \chi =\sup\{\chi (X\cap E):E~\text{~is a Zariski connected irreducible simple
% subgroup of~}~\BG\}.
%$$
%Since there are only finitely many isomorphism classes of
%connected simple subgroups of $\BG$ over the algebraic closure of
%$\BK$, it follows from Theorem \ref{thm:Bezout} that $\chi$ is
%finite. Thus, by Lemma \ref{Bezout} there is some constant
%$N=N(\chi)$ such that whenever $\gL\subset\BG$ is a set containing
%$1$ which generates a subgroup whose Zariski closure is simple and
%irreducible and not contained in $X$, then $\gL^N$ is not
%contained in $X$. It follows that for some word $W$ of length at
%most $N$ in $\Sigma^g$, we have $\| u^W-1\|>\epsilon_G$, for
%otherwise $\langle u^W:W\text{~is a word of length~}\leq
%N\text{~in~}\Sigma^g\rangle$ would be a unipotent group and hence
%some conjugate of it would lie in $U$, and since the corresponding
%conjugate of $\Sigma^g$ generates a group whose Zariski closure
%$\BF^g$ is simple and irreducible, this contradicts the property
%of $N(\chi)$. It follows that equation (\ref{estimate pi(g)})
%holds with $N=N(\chi)$.

We now give an alternative proof which holds in arbitrary
characteristic. Let $U$ be a maximal unipotent subgroup of $\BG$.
For any $u\in U\setminus\{1\}$ let $Y_u=\{ h\in\BG:u^h\in U\}$.
Then $Y_u$ is a proper algebraic subset of $\BG$, and one easily
sees that $\chi (Y_u)$ is bounded independently of $u$, by some $\chi$ say. Now if $\BE$
is an irreducible subgroup of $\BG$, i.e. not contained in a proper
parabolic subgroup, then $\cap_{h\in \BE}U^h$ is trivial. It follows
that $\BE\varsubsetneq Y_u$ for any $u\in U\setminus\{1\}$.
% Let
%$$
% \chi =\sup\{\chi (Y_u\cap \BE):\BE~\text{~is a semisimple algebraic
% subgroup of~}~\BG, u\in U\setminus\{1\} \}.
%$$
%$Since there are only finitely many isomorphism classes of
%$semisimple algebraic subgroups of $\BG$, it follows from Theorem \ref{thm:Bezout} that $\chi$ is finite.
Thus Lemma \ref{Bezout} yields a constant
$N=N(\chi)$ such that
% whenever $\gL\subset\BG$ is a set containing
%$1$ which generates a subgroup whose Zariski closure is semisimple and irreducible then $\gL^N$ is not
% contained in $Y_u$ for any $u\in
%U\setminus\{1\}$.
%It follows that for
some word $W$ of length at
most $N$ in $\Sigma^g$ satisfies $\| u^W-1\|>\epsilon_G$, where $u$ is the element in $(\ref{eq18})$. For
otherwise, $\langle u^W:W\text{~is a word of length~}\leq
N\text{~in~}\Sigma^g\rangle$ would be a unipotent group and hence
some conjugate of it would lie in $U$, and since the corresponding
conjugate of $\Sigma^g$ generates a group whose Zariski closure
$\BE$ is irreducible in $\BG$, this contradicts the property of
$N(\chi)$. It follows that equation (\ref{estimate pi(g)}) holds
with $N=N(\chi)$.
\end{proof}

%------------------------------------------------------------------

\subsection{An $S$--arithmetic version of the Comparison Lemma}\label{subMini}
Let $\BK,S,\BG,G,\gC,d$ be as in the beginning of this Section (we
do not assume that $\BG$ is isotropic over $\BK$). The goal of the
remaining part of this section is to prove the following
arithmetic version of Lemma \ref{prop:constant-c}.

\begin{prop}
\label{prop:gamma}\label{thm:gamma} For some constant $r$,
depending only on $\BG$, $\BK$, and $S$,
%$\Gamma={\mathbb{G}}({\mathcal{O}}_\BK(S))$
we
have that for any finite subset $\Sigma \subset \Gamma={\mathbb{G}}({\mathcal{O}}_\BK(S)) $ ($\Sigma
\ni id$) generating a subgroup whose Zariski closure $\BF$ is
irreducible in $\BG$, there is an element $\gamma \in
\Gamma $ such that $\Vert \Sigma ^{\gamma }\Vert \leq \Lambda
(\Sigma ^{r})$.
\end{prop}

\begin{rem}
In the proof of Theorem \ref{thm2} in the next section we will
apply Proposition \ref{thm:gamma} only in the case where
$\BG=\SL_d$ and $\gC=\SL_d(\OO_\BK(S))$.
\end{rem}

%We will denote by $F$ (resp. $\BF$) the Zariski closure of
%$\langle\gS\rangle$ in $G$ (resp. in $\BG$).

Note that if $\alpha\in{\mathbb{G}}({\mathcal{O}}_\BK (S))$,
$\Lambda (\alpha )=1$ if and only if all the eigenvalues of
$\alpha$ are roots of unity, i.e. if and only if $\alpha$ has
finite order. Moreover, there is a positive constant $\tau>1$ such
that if $\alpha\in{\mathbb{G}}({\mathcal{O}}_\BK (S))$ has
$\Lambda (\alpha )>1$
then $\Lambda (\alpha )\geq\tau$. This follows from the fact that ${\mathcal{O%
}}_\BK (S)$ embeds discretely in $\prod_{v\in S}{\mathbb{K}}_v$.
Moreover, the Zariski closure $Y$ of the set of torsion elements
in $\Gamma$ is a proper algebraic subvariety of $\BG$ (there is an
upper bound of the order of torsion elements in $\Gamma$, see
Proposition 2.5. \cite{Tits}).
%Since there are only finitely many
%isomorphism classes of semisimple algebraic subgroups of
%$\BG$ it follows
%from Theorem \ref{thm:Bezout} that the number $\chi (\BF\cap Y)$ is unifomrly bounded.
Hence Lemma \ref{Bezout} implies that $\Sigma^{\overline{r}}$
contains a non torsion element, where $\overline{r}$ is some
integer independent of $\gS$.
%, and this
%$\overline{r}$ is the same for all choices of a finite subset
%$\Sigma \ni id$ in $\Gamma$ generating a Zariski-dense subgroup.
This shows that $\Lambda(\Sigma^{n})\geq \tau$ for all
$n\geq\overline{r}$. We can therefore reformulate the Comparison
Lemma (\ref{prop:constant-c}) as follows, omitting the
multiplicative constant.

\begin{lem}
\label{r1} For some constant $r'$, depending only on $\BG$, $\BK$, and $S$,
%$\Gamma={\mathbb{G}}({\mathcal{O}}_\BK(S))$
we
have that for any finite subset $\Sigma \subset \Gamma={\mathbb{G}}({\mathcal{O}}_\BK(S)) $ ($\Sigma
\ni id$) generating a subgroup whose Zariski closure $\BF$ is
irreducible in $\BG$, there is an element $h\in H$ such that
$\Vert \Sigma ^{h}\Vert \leq \Lambda (\Sigma ^{r^{\prime }})$.
\end{lem}

In order to derive Proposition \ref{prop:gamma} from Lemma
\ref{r1} we will first replace the conjugating element $h\in H$ by
an element $g\in G$ (of course this step is unnecessary when $G=H$
which is the situation in the proof of Theorem \ref{thm2}). The
second part of the proof which consists in replacing $g$ by some
$\gamma\in\Gamma$ relies on Theorem \ref{thm:DG}.
%The first step follows from the
%fact that the projection on a convex subset in a $CAT(0)$ space is
%distance non-increasing (i.e. $1$-Lipschitz). The second relies on
%the finiteness of the volume of $G/\gC$, namely the Borel
%Harish-Chandra theorem and an effective version of the
%Kazhdan--Margulis theorem.

\medskip

%---------------------------------

\subsubsection{Step 1: Projection to a homogeneous subspace}

By Theorem \ref{thm:ML} we may identify the symmetric space (resp. affine building) of
$G_v$ with a convex subset $C$ of $X_v$ of the form $G_v\cdot x_1$
(resp. $G_v\cdot\gs_1$) for some point $x_1$ (resp. some cell $\gs_1 \ni x_1$) in $X_v$.

Since $X_v$ is a CAT(0) space, the projection to the nearest point
$P_C:X_v\to C$ is $1$--Lipschitz. Let
$h\in\text{SL}_d({\mathbb{K}}_v)$ be the element from Lemma
\ref{r1}, let $x=P_C(h\cdot x_0)$ and let $g_v\in G_v$ be an
element such that $g_v\cdot x=x_1$ (resp. $g_v\cdot x\in \gs_1$). In
any case, we have $d(x_1,g_v\cdot x)\leq 1$. Since $\gS\subset
G_v$, it preserves $C$ and since $P_C$ is $1$--Lipschitz we have
$d_\gS(x)\leq d_\gS(h\cdot x_0)$, where
$d_\gS(x)=\max_{\gc\in\gS}d(x,\gc\cdot x)$. We get $d_\gS (g_v^{-1}\cdot x_1)\leq
d_\gS(h\cdot x_0)+2$, and finally we obtain:
$$
 d_\gS (g_v^{-1}\cdot x_0)\leq
 d_\gS(h\cdot x_0)+2+2d(x_0,x_1).
$$
With Lemma \ref{lem:comp} we can translate this to:
$\|\gS^{g_v}\|\leq\|\gS^h\|^b$ for some constant $b>0$. Repeating
this argument for every $v \in S$, we get from Lemma \ref{r1}:

\begin{cor}
\label{r2} For some constant $r^{\prime \prime }$ (independent of $\gS$) we have
\begin{equation*}
\Vert \Sigma ^{g}\Vert \leq \Lambda (\Sigma ^{r^{\prime \prime
}}).
\end{equation*}
\end{cor}

%-----------------------------------------------------------------------------

\subsubsection{Step 2: Finding a relatively close point in a given
$\gC$--orbit} We will now explain how to replace $g=(g_v)\in G$
by some $\gc\in\gC$ and obtain the proof of Proposition
\ref{prop:gamma}.

Assume first that ${\mathbb{G}}$ is ${\mathbb{K}}$--anisotropic, i.e. that $%
G/\Gamma$ is compact. Let $\Omega$ be a fixed bounded fundamental
domain for
$\Gamma$ in $G$ and let $\gamma\in \Gamma$ be the unique element such that $%
g\in \Omega\gc$. Write
\begin{equation*}
c=\max \{\| f\| :f\in \Omega\cup\Omega^{-1}\},
\end{equation*}
then Theorem \ref{prop:gamma} holds with $r=r''(1+2\log_{\tau}c)$.

Next assume that $\BG$ is $\BK$--isotropic. By equation
(\ref{estimate pi(g)})
\begin{equation*}
\|\pi(g)\|\leq \big( l_\gC\frac{\|\Sigma^g\|^{2N}}{\epsilon_G}\big)%
^{1/k_{\gC}}\leq \big( l_\gC\frac{\Lambda (\Sigma^{r''})^{2N}}{\epsilon_G}\big)%
^{1/k_{\gC}}.
\end{equation*}
Which means that for some $\gamma\in\Gamma$
\begin{equation*}
\|g\gamma^{-1}\|\leq \big( l_\gC\frac{\Lambda (\Sigma^{r''})^{2N}}{\epsilon_G}%
\big)^{1/k_{\gC}} \leq \Lambda (\Sigma^{r''} )^{\frac{2N}{k_{\gC}}+\frac{1}{k_{\gC}}\log_\gt\frac{%
l_\gC}{\epsilon_G}}.
\end{equation*}
and therefore
\begin{equation*}
\|\Sigma^\gc\|=\|\Sigma^{\gamma
g^{-1}g}\|\leq\|\Sigma^g\|\|g\gamma^{-1}\|\|\gamma g^{-1}\|\leq \Lambda
(\Sigma^r)
\end{equation*}
for some computable constant $r$. \qed

%-----------------------------------------------------------------------

\section{Construction of the ping-pong players}\label{sec:ping-pong}

In this section we will construct two bounded words in the
alphabet $\Sigma$ that will play ping-pong on some projective
space and hence will be independent. This will prove Theorem
\ref{thm2}. Since the detailed proof below is somewhat technical
we refer the reader to \cite{note} for an outline of the main
ideas.

Let $\gS\subset\SL_d(\OO_\BK(S))$ be as in the statement of
Theorem \ref{thm2}. Inconsistently with the previous section we
will denote by $\BG$ the Zariski closure of $\langle\gS\rangle$ in
$\SL_d$, and $\gC=\BG (\OO_\BK (S))$. Then $\BG$ is a semisimple
irreducible subgroup in $\SL_d$ and $\gC$ is an $S$--arithmetic
subgroup of $\BG$. We let $G=\prod_{v\in S}\BG(K_v)$ and identify
$\gC$ via the diagonal embedding with the corresponding
$S$--arithmetic lattice in $G$.

In this section, whenever we say that some quantity is a \textit{constant}, we mean that it
may depend only on $d$, $\BK$ and $S$.

The following proposition will allow us to assume that $\gS$ is
finite, hence compact. Let $s\in\BN$ be a constant. We will
specify some condition on $s$ in Paragraph \ref{sbs:step4} (Step
(4)), for the moment we only require it to be at least $r$, the
constant from Proposition \ref{prop:gamma}.

\begin{prop}\label{sigma-two}
There is a constant $f$, such that for
any subset $1\in\gS\subset\gC$ which generates a Zariski dense
subgroup of $\BG$, there is a subset $\gS'$ of $\gS^f$ of
cardinality $\dim (\BG )$ such that:
\begin{enumerate}
\item The Zariski closure $\overline{\langle\gS'\rangle}^Z$ of the group generated by $\gS'$ equals $\BG$, and

\item $(\gS')^s$ consists of semisimple elements.
\end{enumerate}
\end{prop}

Recall the following fact:

\begin{lem}\label{Borel}\label{lem:dominant} \textbf{(see Borel \cite{Bor})} Let $\BG$ be a connected semisimple algebraic
group, and $k\geq 2$ an integer. If $W$ is a non-trivial word in
the free group $F_k$, then the corresponding map $W: \BG^k
\rightarrow \BG$ is dominant.
\end{lem}

\begin{proof}[Proof of Proposition \ref{sigma-two}]
Let $k=\dim (\BG)$, let $W_1,\ldots,W_t$ be all the reduced words
in $F_k$ of length $\leq s$, and consider the map
$w:\BG^k\to\BG^t$ defined by substitution in $(W_1,\ldots,W_t)$.
Let $\Phi\subset\BG$ be a Zariski open subset which consists of
semisimple elements. We shall construct inductively elements
$\gs_i,~i=1,\ldots$ in a bounded power of $\gS$ which $\forall i$
satisfy:
\begin{itemize}
\item There are some $g_{i+1},\ldots,g_k\in\BG$ such that
$w(\gs_1,\ldots,\gs_i,g_{i+1},\ldots,g_k)\in \Phi^t$.

\item $\dim (\overline{\langle\gs_1,\ldots\gs_i\rangle}^Z)\geq i$.
\end{itemize}

In order to construct $\gs_1$ choose some $(g_1,g_2,\ldots,g_k)\in
w^{-1}(\Phi^t)$, which is non-empty by Lemma \ref{lem:dominant} ,
and define
$$
 V_1=\{g\in\BG:w(g,g_2,\ldots,g_k)\in\BG^t\setminus\Phi^t\}.
$$
As noted before Lemma \ref{r1}, the Zariski closure $X$ of the
elements in $\gC$ whose projection to one of the factors of $\BG$
is torsion is a proper subvariety of $\BG$. Let $N_1$ be the
constant obtained from Lemma \ref{Bezout} applied to $X\cup V_1$
and take $\gs_1\in\gS^{N_1}\setminus (X\cup V_1)$. It is
straightforward to check that $\chi(V_1)$ (i.e. the sum of the
degrees and dimensions of the irreducible components of $V_1$) can
be bounded independently of the choice of $(g_2,\ldots,g_k)$ and
hence that $N_1$ can be taken to be a constant. Finally, since
$\gs_1$ has infinite order $\overline{\langle\gs_1\rangle}^Z$ has
positive dimension.

To explain the $i$'th step let us suppose that
$\gs_1,\ldots,\gs_{i-1}$ were already constructed. Since
$\gs_1,\ldots,\gs_{i-1}$ are assumed to satisfy the requirements
above, we can chose some new $g_{i+1},\ldots,g_k\in \BG$ for which
the algebraic set
$$
 V_i=\{g\in\BG:w(\gs_1\ldots,\gs_{i-1},g,g_{i+1},\ldots,g_k)\in\BG^t\setminus\Phi^t\}
$$
is proper. Additionally, the Zariski connected group
$\BG_i=(\overline{\langle\gs_1,\ldots,\gs_{i-1}\rangle}^Z)^\circ$
cannot be proper normal since by the properties of $\gs_1$ it
projects non-trivially to each simple factor of the semisimple
group $\BG$. If $\BG_i=\BG$ take $\gs_i=1$ and otherwise take $\gd_i\in \gS\setminus N_\BG(\BG_i)$, let $N_i$ be the constant obtained from Lemma \ref{Bezout} applied to $V_i\cup \gd V_i$, chose $\gs_i'\in\gS^{N_i}\setminus
(V_i\cup\gd_i V_i)$, and set $\gs_i=\gs_i'$ if $\gs_i'\notin N_\BG(\BG_i)$ and $\gs_i=\gd_i^{-1}\gs_i'$ otherwise.
Again $N_i$ can be taken to be a constant (independent of the previous choice of $\gs_j,~j<i$, the
choice of $g_j,~j>i$ and the choice of $\gd_i$, since $\chi (V_i\cup\gd_iV_i)$ too can be bounded by a constant). Finally, since $\gs_i$ does not normalize $\BG_i$, $\dim(\overline{\langle\gs_1,\ldots,\gs_i\rangle}^Z>\dim(\overline{\langle\gs_1,\ldots,\gs_{i-1}\rangle}^Z$.
\end{proof}

We will therefore assume that $\gS$ itself is finite and $\gS^s$
consists of semisimple elements (where $s\geq r)$. Applying Proposition
\ref{prop:gamma} we see that up to changing $\gS$ into $\gS^{\gc}$
for some $\gamma\in\Gamma$, we may assume that $\Lambda (A_0)\geq
\|\gS\|$ for some $A_0\in\gS^r$.

We will now fix once and for all a place $v\in S$ for which
$\Lambda_{v}(A_0)=\Lambda (A_0)$. The local field
${\mathbb{K}}_{v}$ has only finitely many extensions of degree at
most $d!$. Let $\tilde{{\mathbb{K}}}_{v}$ be their compositum,
then any semisimple element in $\text{SL}_d({\mathbb{K}}_{v})$ is
diagonalizable in $\text{SL}_d(\tilde{{\mathbb{K}}}_{v})$. Similarly, let $%
\tilde{{\mathbb{K}}}$ be the splitting field of $A_0$, and let
$\tilde{S}$ be the set of all places of $\tilde{{\mathbb{K}}}$
extending elements of $S$.

By passing to a suitable wedge power representation $V=\gL^i\BK^d$
for some $i$, $1 \leq i \leq d-1$,
%\footnote{In the $0$ characteristic case we
%can actually take an irreducible sub-representation because $\BG$
%is semisimple and hence any representation is completely reducible.}),
we may assume that $A_0$ has a unique eigenvalue $\alpha_1(A)$ of
maximal $v$--absolute value  and that the ratio between
$\alpha_1(A_0)$ and the second largest eigenvalue $\alpha_2(A_0)$
satisfies
\begin{equation*}
\Lambda(A_0)^{d}\geq\Big|
\frac{\alpha_1(A_0)}{\alpha_2(A_0)}\Big|_{v}\geq \Lambda
(A_0)^{\frac{1}{d}} \geq \tau^{1/d},
\end{equation*}
where $\tau$ is the constant introduced in the proof of
Proposition \ref{prop:gamma}. Note that the norm of a matrix in a
wedge power representation such as $V$ is bounded by its original
norm to the power $d$. Thus, we have
\begin{equation}  \label{lambda1/lambda2}
|{\alpha_1(A_0)}/{\alpha_2(A_0)}|_{v}^{d^2}\geq \| \gS
\|_{End(V)}.
\end{equation}

We will set $n=dimV$ the dimension of the new representation. Note
that $n \leq 2^d$. Note also that in the canonical basis of the
wedge power space, the matrix elements from $\gS$ (viewed as
matrices in $\SL_n(\BK)$) are still in $\OO_\BK (S)$. Finally
observe that $V$ may not be $\BG$--irreducible. This is not a
fundamental problem. However to keep exposition as simple as
possible we will assume throughout that $V=\BK^{n}_v$ is an
irreducible $\BG$--space with $A_0$ and $\gS$ with matrix
coefficients in $\OO_\BK (S)$ and satisfying the two inequalities
above. At the end we will indicate the changes to be made to
accomodate with the fact that $\gL^i\BK^d$ is not irreducible in
general.

Working with the corresponding projective representation over
$\tilde{{\mathbb{K}}}_{v}$ we will now produce two ping--pong
players in four steps. In the first we will construct a proximal
element, in the second a very contracting one and in the third a
very proximal one. Then we will find a suitable conjugate of the
very proximal element and obtain in this way a second ping--pong
partner.

%----------------------------------------------
\subsection{Step 1}

We set $r_0=rd^2$. Let $\{\hat{u}_i\}$ be a basis of
$\tilde{{\mathbb{K}}}_v^{n}$ consisting of
normalized eigenvectors of $A_0$ with corresponding eigenvalues $\{\alpha_i\}$%
, such that whenever $\alpha_i=\alpha_j$ the vectors $\hat{u}_i$ and $\hat{u}%
_j$ are orthogonal\footnote{In the non-Archimedean case this is
simply taken to mean that $\|\hat u_i-\hat u_j\|=1$.}, and let
$\hat u_i^\perp$ denote the hyperplane spanned by
$\{\hat{u}_j:j\neq i\}$.

\begin{lem}
\label{const:r1} For some constant $r_{1}\in {\mathbb{N}}$,
depending only on $\Gamma $,
\begin{equation*}
d(\hat{u_{i}},\hat{u}_{i}^{\perp })\geq |\frac{\alpha _{1}}{\alpha _{2}}%
|_{v}^{-r_{1}}
\end{equation*}
for $i=1,\ldots n$.
\end{lem}

\begin{proof}
First note that since $|\ga_i-\ga_j|_w\leq 2\Lambda (A_0)$ for any
$w\in\ti S$ and $|\ga_i-\ga_j|_w\leq 1$ for any $w\notin \ti S$,
it follows from the product formula that if $\ga_i\neq\ga_j$ then
$$
 |\ga_i-\ga_j|\geq (2\Lambda (A_0))^{-|\ti S|}\geq \Lambda
 (A_0)^{-|\ti S|(1+\log_{\gt}2)}\geq |\frac{\ga_1}{\ga_2}|_v^{-d|\ti S|(1+\log_{\gt}2)}=
 |\frac{\ga_1}{\ga_2}|_v^{-t_0}
$$
where $t_0=d|\ti S|(1+\log_{\gt}2)$. Note also that $|\ti S|\leq
d!|S|$.

Next, observe that it is enough to show that for some constant
$r_1'$,
\begin{equation*}
 d(\overrightarrow{u_i},\text{span}\{ \hat{u}_j:\ga_j\neq\ga_i\})\geq |\frac{\ga_1}{\ga_2}|_v^{-r_1'}
\end{equation*}
for any $i$ and any unit vector
$\overrightarrow{u}_i\in\text{span}\{\hat u_j:\ga_j=\ga_i\}$. This
in turn will follow from the next claim which we will prove by
induction on $k$:

\medskip

\noindent {\bf Claim.} For any $k$ there is a positive constant
$t_k$ such that if
$\overrightarrow{u}\in\text{span}\{\hat{u_j}:j\in
I,\ga_j\neq\ga_i\}$ where $I$ is a set of indices with $\dim
(\text{span}\{\hat{u_j}:j\in I,\ga_j\neq\ga_i\})=k$ then
$\|\overrightarrow{u}_i-\overrightarrow{u}\|\geq|\frac{\ga_1}{\ga_2}|_v^{-t_k}$
for any unit vector $\overrightarrow{u}_i\in\text{span}\{\hat
u_j:\ga_j=\ga_i\}$.

\medskip

For $k=1$ we can write $\overrightarrow{u}=\lambda\hat{u}_j$, so
$$
 A_0(\overrightarrow{u}_i-\lambda\hat{u}_j)=(\ga_i-\ga_j)\overrightarrow{u}_i+
 \ga_j(\overrightarrow{u}_i-\lambda\hat u_j)
$$
i.e.
$$
 (A_0-\ga_j)(\overrightarrow{u}_i-\lambda\hat{u}_j)=
 (\ga_i-\ga_j)\overrightarrow{u}_i,
$$
which implies that (recall $r_0=rd^2$)
$$
 \|\overrightarrow{u}_i-\lambda\hat
 u_j\|_v\geq\frac{|\ga_i-\ga_j|_v}{\|A_0\|_v+|\ga_j|_v}\geq
 |\frac{\ga_1}{\ga_2}|_v^{-t_0-(r_0+d\log_{\gt}2)}:=
 |\frac{\ga_1}{\ga_2}|_v^{-t_1}.
$$

Now suppose $k>1$. We can write
$\overrightarrow{u}=\sum\lambda_j\overrightarrow{u}_j$ where the
$\overrightarrow{u}_j$'s are normalized eigenvectors of different
eigenvalues. Abusing indices, we will assume that
$\overrightarrow{u}_j$ corresponds to the eigenvalue $\ga_j$. Now
$$
 A_0(\overrightarrow{u}_i-\sum\lambda_j\overrightarrow{u}_j)=
 \ga_i(\overrightarrow{u}_i-\sum\lambda_j\overrightarrow{u}_j)+
 \sum_j(\ga_i-\ga_j)\lambda_j\overrightarrow{u}_j,
$$
therefore
$$
 (A_0-\ga_i)(\overrightarrow{u}_i-\sum\lambda_j\overrightarrow{u}_j)=
 \sum_j(\ga_i-\ga_j)\lambda_j\overrightarrow{u}_j.
$$
Note that we may assume that $\|\overrightarrow{u}\|_v\geq 1/2$,
for otherwise the statement is obvious, and hence for some $j_0$,
$|\lambda_{j_0}|_v\geq 1/(2n)$ and by the induction hypothesis
$$
 \|\sum_j(\ga_i-\ga_j)\lambda_j\overrightarrow{u}_j\|\geq
 |\lambda_{j_0}|_v|\frac{\ga_1}{\ga_2}|_v^{-t_0-t_{k-1}}\geq
 \frac{1}{2n}|\frac{\ga_1}{\ga_2}|_v^{-t_0-t_{k-1}}.
$$
It follows that
$$
 \|\overrightarrow{u}_i-\sum\lambda_j\overrightarrow{u}_j\|_v\geq
 \frac{1}{2n}|\frac{\ga_1}{\ga_2}|_v^{-t_0-t_{k-1}}\frac{1}{\|A_0\|_v+|\ga_i|_v}\geq
 |\frac{\ga_1}{\ga_2}|_v^{d\log_{\gt}\frac{1}{2n}-t_0-t_{k-1}
 -(r_0+d\log_{\gt}2)}:=|\frac{\ga_1}{\ga_2}|_v^{-t_k}.
$$
\end{proof}

As a consequence we obtain that for some constant $r_2$, depending
only on $\Gamma$, which we may take $\geq r_1$, we have:

\begin{cor}
\label{D} There is a matrix $D\in
\text{SL}_{n}(\tilde{{\mathbb{K}}}_{v})$ such that:

\begin{itemize}
\item  $\Vert D\Vert^2 ,\Vert D^{-1}\Vert^2 \leq |\frac{\alpha _{1}}{\alpha _{2}}%
|_{v}^{r_{2}}$, and

\item  $A_{0}^{D}=DA_0D^{-1}$ is diagonal.
\end{itemize}
\end{cor}

\begin{proof}
Let $D$ be the matrix defined by the condition $D(\hat
u_i)=e_i,~i=1,\ldots,n$. Clearly $|\text{det}(D^{-1})|_v\leq 1$.
Since $D^{-1}=\text{det}(D^{-1})\text{Adj}(D)$ and since
$\|\text{Adj}(D)\|\leq n! \| D\|^{n-1}$ it is enough to prove that
$\| D\|\leq |\frac{\ga_1}{\ga_2}|_v^{r_2'}$.

Let $\hat{u}$ be a unit vector, and write $\hat u=\sum\lambda_i
e_i$. Then for some $i_0$ we have $|\lambda_{i_0}|_v\geq 1/n$.
Since $D^{-1}(\hat u)=\sum\lambda_i\hat{u}_i$, it follows from the
previous lemma that
$$
 \| D^{-1}(\hat u)\| =\|\sum\lambda_i\hat{u}_i\|=
 \|\lambda_{i_0}\hat{u}_{i_0}+\sum_{j\neq i_0}\lambda_j\hat u_j\|\geq
 \frac{1}{n}|\frac{\ga_1}{\ga_2}|_v^{-r_1}\geq
 |\frac{\ga_1}{\ga_2}|_v^{-r_2'},
$$
i.e. $\| D\|\leq |\frac{\ga_1}{\ga_2}|_v^{r_2'}.$
\end{proof}

We derive the following proposition and thus conclude the first
step in our construction of ping--pong players:

\begin{prop}[The proximal element $A_1$]
Whenever $r_{3}\geq 8r_{2}$, the element $A_{1}=A_{0}^{r_{3}}$ is $(|%
\frac{\alpha _{1}}{\alpha _{2}}|_{v}^{-r_{1}},|\frac{\alpha _{1}}{\alpha _{2}%
}|_{v}^{-(r_{3}/2-2r_{2})})$--proximal with attracting point
$[\hat{u}_{1}]$
and repelling hyperplane $[\hat{u}_{1}^{\perp }]=[{\text{span}(\hat{u}%
_{2},\ldots ,\hat{u}_{n})}]$.
\end{prop}

\begin{proof}
The diagonal matrix $DA_0^{r_3}D^{-1}$ is obviously
$|\frac{\ga_1}{\ga_2}|_v^{-r_3/2}$--contracting with attracting
point $[{e}_1]$ and repelling hyperplane $[{\text{span}(e_2,\ldots
,e_n)}]$. Since $\| D\|,\| D^{-1}\|\leq
|\frac{\ga_1}{\ga_2}|_v^{r_2}$, $D$ is
$|\frac{\ga_1}{\ga_2}|_v^{2r_2}$ bi-Lipschitz. It follows that
$A_0^{r_3}$ is
$|\frac{\ga_1}{\ga_2}|_v^{-(r_3/2-2r_2)}$--contracting. Finally,
Lemma \ref{const:r1} implies that $d([\hat{u}_1],[\hat
u_1^{\perp}])\geq |\frac{\ga_1}{\ga_2}|_v^{-r_1}$
\end{proof}

%-------------------------------------------------------

\subsection{Step 2}

Our next goal is to build a very contracting element out of the
matrix $A_1$. To achieve this, we will find some bounded
word $B_1$ in $\gS$ which will be in ``general position"
with respect to $A_1$. Then $A_2=A_1^{r_7}B_1A_1^{-r_7}$ will be
our candidate. In this process we will ``lose" the information we
have on the position of the repelling neighborhoods. However we
will still have a good control on the positions of the attracting
points of $A_2$ and $A_2^{-1}$, a control which will turn crucial
in the following step when producing a very proximal element
$A_3$. The key idea is that while $B_1$ sends the eigen-directions
of $A_1$ away from the eigen-hyperplanes of $A_1$, we can estimate
this quantitatively by giving an explicit lower bound. In order to
formulate a precise statement, we will need to introduce another
basis of eigenvectors for $A_1$.

\begin{lem}
\label{nice-vectors} For each $k\leq n$ there is an eigenvector $%
\overrightarrow{u}_{k}\in \overline{{\mathbb{K}}}^{n}$ for $A_{0}$
with corresponding eigenvalue $\alpha _{k}$ whose coordinates are
$\ti S$--integers and whose $w$--norm is at most $|\alpha
_{1}/\alpha _{2}|_{v}^{r_{4}}$ for any $w\in \ti S$, where $r_{4}$
is some constant depending only on $r_{0},d$ and the size of $S$.
\end{lem}

\begin{proof}
Recall from inequality (\ref{lambda1/lambda2}) that for each $w\in
S$ we have $\| A_0\|_w\leq |\ga_1/\ga_2|_v^{r_0}$ (where
$r_0=rd^2$). Suppose that $\ga_i$ has multiplicity $k$, say
$\ga_i=\ga_{i+1}=\ldots =\ga_{i+k-1}$, then we can pick $k$
indices between $1$ and $n$ such that the $(n-k)\times (n-k)$
matrix obtained by restricting $A_0 -\ga_i$ to the remaining
indices is invertible. We can then define
$\overrightarrow{u_{i+j}},~j\leq k-1$ to be the eigenvector of
$\ga_{i+j}=\ga_i$ whose entries corresponding to the chosen $k$
indices are all $0$ except the $(j-1)$'th one which equals the
determinant of the $(n-k)\times (n-k)$ submatrix. Solving the
corresponding linear equation, it is easy to verify that these
vectors satisfy the requirement with respect to some bounded
constant $r_4$.
\end{proof}

In analogy to our previous notations, we will denote by $\overrightarrow{u_i%
}^\bot$ the span of the $\overrightarrow{u_j}$'s, $j\neq i$. Note
that since $\alpha_1$ has multiplicity one, we have $[\hat
u_1]=[\overrightarrow{u_1}]$ and $[\hat
u_1^\bot]=[\overrightarrow{u_1}^\bot]$.

\begin{defn}\label{defn:general-position}
Let $N$ be an integer and $v_1,\ldots,v_n\in\overline{\BK}^n$ a
basis. We will say that a matrix $C\in\SL_n(\overline{\BK} )$ is
in $N$--{\it general position} with respect to $\{
v_1,\ldots,v_n\}$ if
\begin{itemize}
\item for any $1\leq i,j\leq n$, not necessarily distinct, both
vectors $Cv_i$ and $C^{-1}v_i$ do not lie in the hyperplane
spanned by $\{v_k\}_{k\neq j}$, and

\item for any $n$ integers $1\leq i_1<\ldots <i_n\leq N$ and any
$1\leq j\leq n$ the vectors $C^{i_1}v_j,\ldots,C^{i_n}v_j$ are
linearly independent.
\end{itemize}
\end{defn}

For a fixed $N$, the varieties
$$
 X(N,v_1,\ldots,v_n)=\{g\in\SL_n(\overline{\BK} ):g~\text{is {\it not} in $N$--general position w.r.t.}~\{v_i\}_{i=1}^n\}
$$
are all conjugate inside $\SL_n(\overline{\BK} )$.
%It follows from Bezout's
%theorem that there is an upper bound on $\chi
%(X(N,v_1,\ldots,v_n)\cap\BG)$ which is independent of the $v_i$.
Since $\BG$ is Zariski connected and irreducible, one can derive
that $X(N,v_1,\ldots,v_n)\cap\BG$ is a proper subvariety of $\BG$.
Hence by Lemma \ref{Bezout} for any $N$ there is a constant
$m_2(N)$ such that for any set $\Omega$ which generates a Zariski
dense subgroup of $\BG$, and any basis $\{ v_i\}_{i=1}^n$ of
$\BK^n$, there is an element in $\Omega^{m_2(N)}$ which is in
$N$--general position with respect to $\{v_i\}_{i=1}^n$. In
particular we may find $B_1\in \gS^{m_2}$ (with $m_2=m_2(2n-1)$)
which is in $(2n-1)$--general position with respect to
$\{\overrightarrow{u}_i\}_{i=1}^n$.

In the proof of Proposition \ref{very-contracting} we will make
use of the following lemma only for $i=n$ and $j=1$.

\begin{lem}
\label{r7}\label{const:r5} For some positive bounded constant
$r_{5}$ we have
\begin{equation*}
d((B_1^{\pm 1})\cdot [\overrightarrow{u_{i}}],[\overrightarrow{u_{j}^{\perp }}%
])>|\frac{\alpha _{1}}{\alpha _{2}}|_{v}^{-r_{5}},
\end{equation*}
for any $i,j\leq n$.
\end{lem}

\begin{proof}
For each $w\in\ti S$, the $w$--absolute values of the coordinates
of $B_1(\overrightarrow{u}_i)$ are at most
$|\ga_1/\ga_2|_v^{m_2r_0+r_4}$. Consider the determinant
$$
 \mathcal{D}_{\pm 1}=\text{det}(B^{\pm 1}(\overrightarrow{u}_i),\overrightarrow{u}_1,\ldots,
 \overrightarrow{u}_{j-1},\overrightarrow{u}_{j+1},\ldots,\overrightarrow{u}_n).
$$
This is again an $\ti S$--integer and its $w$--absolute value is
at most $|\ga_1/\ga_2|_v^{m_2r_0+nr_4}$. Since $B_1$ is in general
position with respect to $\{\overrightarrow{u}_i\}_{i=1}^n$ we
have $\mathcal{D}_{\pm 1}\neq 0$. By the product formula
$\prod_{\text{all places}}|\mathcal{D}_{\pm 1}|_w=1$ and hence
$\prod_{w\in\ti S}|\mathcal{D}_{\pm 1}|_w\geq 1$. It follows that
$$
 |\mathcal{D}_{\pm 1}|_v\geq |\ga_1/\ga_2|_v^{-(m_2r_0+nr_4)|\ti S|}.
$$
Now since all the vectors involved in this determinant have
$v$--norm at most $|\ga_1/\ga_2|_v^{m_2r_0+r_4}$, the distance
between each of them to the hyperplane spanned by the others is at
least
$$
 \frac{|\mathcal{D}_{\pm 1}|_v}{|\ga_1/\ga_2|_v^{(m_2r_0+r_4)(n-1)}}\geq
 |\ga_1/\ga_2|_v^{-(m_2r_0+nr_4)(|\ti S|+n-1)}.
$$
The lemma follows.
\end{proof}

We will also need the following:

\begin{lem}
There exists some $\epsilon =\epsilon (n)$, such that if $d=\text{diag}%
(d_{1},\dots ,d_{n})\in \text{SL}_{n}(\tilde{{\mathbb{K}}}_{v})$
is a
diagonal matrix with $d_{1}\geq d_{2}\geq \ldots \geq d_{n}$, then $[d]$ is $%
2$--Lipschitz on the $\epsilon$--ball around $[e_{1}]$.
\end{lem}

\begin{proof}
The lemma follows by a direct simple computation. In the
non-Archimedean case a diagonal matrix is $1$--Lipschitz on the
open unit ball around $[e_1]$. In the Archimedean case the same is
true for the metric which is induced on $\BP (\ti\BK_v^n)$ from
the $L^{\infty}$ norm on $\ti\BK_v^n$. Since the renormalization
map from the euclidean unit sphere to the $L^{\infty}$ unit sphere
is $C^1$ around $e_1$ with differential $1$ at $e_1$ it has a
bi-Lipschitz constant arbitrarily close to $1$ in a small
neighborhood of $e_1$. The result follows.
\end{proof}

We are now able to formulate:

\begin{prop}[The very contracting element $A_2$]
\label{very-contracting} For any $r_{6}\in\BN$, there exists
$r_{7}\in\BN$ such that the element
$A_{2}=A_{1}^{r_{7}}B_1A_{1}^{-r_{7}}$ is $|\frac{\alpha
_{1}}{\alpha _{2}}|_{v}^{-r_{6}}$ very contracting, with both
attracting points (of the
element and its inverse) lying in the $|\frac{\alpha _{1}}{\alpha _{2}}%
|_{v}^{-r_{6}}$ ball around $[\hat{u}_{1}]$.
\end{prop}

The proof of Proposition \ref{very-contracting} relies on
Proposition \ref{contracting-properties}, as well as the last two
lemmas:

\begin{proof}[Proof of Proposition \ref{very-contracting}]
Let $r_7\in\BN$ be arbitrary, to be determined later. By the
previous lemma, the diagonal matrix $DA_1^{-r_7}D^{-1}$ is
$2$--Lipschitz on the on the $\gep (n)$--ball around
$[e_n]=D[\hat{u}_n]$. By Corollary \ref{D} $\| D^{\pm 1}\|^2\leq
|\ga_1/\ga_2|_v^{r_2}$ which implies that $D^{\pm 1}$ are
$|\ga_1/\ga_2|_v^{2r_2}$ Lipschitz (on the entire projective
space, see Section \ref{prelim} (iv)). It follows that
$A_1^{-r_7}$ is $2|\ga_1/\ga_2|_v^{4r_2}$ Lipschitz on the
$\gep\cdot |\ga_1/\ga_2|_v^{-2r_2}$ ball around $[\hat u_n]$, or
in other words, that $A_1^{-r_7}$ is
$|\ga_1/\ga_2|_v^{d\log_{\gt}2+4r_2}$ Lipschitz on the
$|\ga_1/\ga_2|_v^{d\log_{\gt}\gep -2r_2}$ ball around $[\hat
u_n]$.

Now since $\| B_1^{\pm 1}\|_v\leq |\ga_1/\ga_2|^{m_2r_0}$, the
matrices $B_1^{\pm 1}$ are $|\ga_1/\ga_2|^{2m_2r_0}$ Lipschitz on
the projective space, and hence the matrices $B_1^{\pm
1}A_1^{-r_7}$ are $|\ga_1/\ga_2|_v^{d\log_{\gt}2+4r_2+2m_2r_0}$
Lipschitz on the $|\ga_1/\ga_2|_v^{d\log_{\gt}\gep -2r_2}$ ball
around $[\hat u_n]$.

Take
$$
 c^*=\max\{ 2r_2-d\log_{\gt}\gep,~~d\log_{\gt}2+4r_2+2m_2r_0+2r_7\},
$$
then the $|\ga_1/\ga_2|_v^{-c^*}$--ball $\gO$ around $[\hat u_n]$
is mapped under $B_1A_1^{-r_7}$ (resp. under $B_1^{-1}A_1^{-r_7}$)
into the $|\ga_1/\ga_2|_v^{-2r_7}$--ball around $B_1[\hat u_n]$
(resp. around $B_1^{-1}[\hat u_n]$). By Lemma \ref{const:r5}
$$
 d(B_1^{\pm 1}[\hat u_n],[\hat u_1^\perp])\geq
 |\ga_1/\ga_2|_v^{-r_5}.
$$
Note that without loss of generality we can set $[\hat u_n]$ to be
equal to $[\overrightarrow{u_n}]$. Also we may assume that
$|\ga_1/\ga_2|_v^{-r_5}<1/\sqrt{2}$ and that $r_7\geq r_5$.
Therefore, $B_1A_1^{-r_7} \gO$ and $B^{-1}_1A_1^{-r_7} \gO$ lie
outside the $|\ga_1/\ga_2|_v^{-2r_5}$ neighborhood of $[\hat
u_1^\bot]$. It follows that both sets $DB_1^{\pm 1}A_1^{-r_7}\gO$
lie outside the $|\ga_1/\ga_2|_v^{-r_5-2r_2}$ neighborhood of
$D[\hat u_1^\bot ]=[\text{span}\{e_2,\ldots, e_n\} ]$. By
Proposition \ref{contracting-properties} $(1)$ applied to the
diagonal matrix $DA_1^{r_7}D^{-1}$, it is
$|\ga_1/\ga_2|_v^{-r_7+2(r_5+2r_2)}$--Lipschitz outside the
$|\ga_1/\ga_2|_v^{-r_5-2r_2}$ neighborhood of
$[\text{span}\{e_2,\ldots, e_n\} ]$, and hence $A_1^{r_7}D^{-1}$
is $|\ga_1/\ga_2|_v^{-r_7+2(r_5+3r_2)}$-Lipschitz there. Thus
$A_1^{r_7}B_1^{\pm 1}A_1^{-r_7}=(A_1^{r_7}D^{-1})D(B_1^{\pm
1}A_1^{-r_7})$ are both
$$
 |\ga_1/\ga_2|_v^{-r_7+c^{**}}-\text{Lipschitz}
$$
on $\gO$, where we have set
$$
 c^{**}=2(r_5+3r_2)+(d\log_{\gt}2+4r_2+2m_2r_0)+2r_2.
$$
It follows from parts (2) and (3) of Proposition
\ref{contracting-properties} that the elements
$A_1^{r_7}B_1^{\pm}A_1^{-r_7}$ are both
$|\ga_1/\ga_2|_v^{\frac{1}{2}[-r_7+c^{**}]}$ contracting. Thus
taking
$$
 r_7\geq 2r_6+c^{**}
$$
we guarantee that $A_1^{r_7}B_1A_1^{-r_7}$ is
$|\ga_1/\ga_2|_v^{-r_6}$ very contracting.

Now suppose further that
$$
 r_7\geq 2\max\{2r_6,c^*\}+c^{**},
$$
then our elements $A_1^{r_7}B_1^{\pm 1}A_1^{-r_7}$ are
$|\ga_1/\ga_2|_v^{-\max\{2r_6,c^*\}}$--very contracting. Moreover
$\gO$ is a $|\ga_1/\ga_2|_v^{-c^*}$--ball, hence contains a point
$p^+$ (resp. a point $p^-$) that is at least
$|\ga_1/\ga_2|_v^{-c^*}$-away from the repelling hyperplane of
$A_1^{r_7}B_1A_1^{-r_7}$  (resp. of
$A_1^{r_7}B_1^{-1}A_1^{-r_7}$). It follows that
$A_1^{r_7}B_1^{\pm}A_1^{-r_7}$ maps the points $p^\pm$
respectively into the $|\ga_1/\ga_2|_v^{-2r_6}$--ball around the
corresponding attracting points $t^\pm$ of $A_1^{r_7}B_1^{\pm
1}A_1^{-r_7}$, i.e.
$$
 d(A_1^{r_7}B_1^{\pm 1}A_1^{-r_7}
 (p^\pm),t^\pm)\leq |\ga_1/\ga_2|_v^{-2r_6}.
$$

Additionally the element $A_1^{r_7}$ is
$|\ga_1/\ga_2|_v^{-r_7/2+4r_2}$--contracting with attracting point
$[\hat u_1]$ and repelling hyperplane $[\hat u_1^\bot]$, and since
the point $B_1[\hat u_n]$ lies outside the
$|\ga_1/\ga_2|_v^{-r_5}$ neighborhood of $[\hat u_1^\bot]$,
assuming further that $r_7\geq 2r_5+8r_2$, we get that this point
is mapped under $A_1^{r_7}$ to the
$|\ga_1/\ga_2|_v^{-r_7/2+4r_2}$--ball around $[\hat u_1]$. We
conclude that $[\hat u_n] \in \gO$ is mapped under
$A_1^{r_7}B_1A_1^{-r_7}$ into the
$|\ga_1/\ga_2|_v^{-r_7/2+4r_2}$--ball around $[\hat u_1]$.

Finally since $A_1^{r_7}B_1A_1^{-r_7}$ is
$|\ga_1/\ga_2|_v^{-r_7+c^{**}}$ Lipschitz on $\gO$, we get that
\begin{eqnarray*}
 d(t^+,[\hat u_1])
 %&\leq&\\
 &\leq& d(t^+,A_1^{r_7}B_1A_1^{-r_7}
 p^+)+\text{diam}(A_1^{r_7}B_1A_1^{-r_7}\gO) +d(A_1^{r_7}B_1A_1^{-r_7}[\hat
 u_n],[\hat u_1])\\
 &\leq& |\ga_1/\ga_2|_v^{-2r_6} +2|\ga_1/\ga_2|_v^{-r_7+c^{**}-c^*}
 +|\ga_1/\ga_2|_v^{-r_7/2+4r_2}.
\end{eqnarray*}
By choosing $r_7$ sufficiently large, we can make the last
quantity smaller the $|\ga_1/\ga_2|_v^{-r_6}$, that is
$$
 d(t^+,[\hat u_1])\leq |\ga_1/\ga_2|_v^{-r_6}.
$$

The same computation with $t^-,p^-$ replacing $t^+,p^+$ gives
$d(t^-,[\hat u_1])\leq |\ga_1/\ga_2|_v^{-r_6}$. This finishes the
proof of the proposition.
\end{proof}

%-------------------------------------------------------------------------------------------------------------

\subsection{Step 3}

Our next step is to use $A_2=A_1^{r_7}B_1A_1^{-r_7}$ to build a
very proximal element. Note that we haven't specified any
condition on the constants $r_6,r_7$ from Lemma \ref
{very-contracting} yet. We will show that for some suitable $k\leq
2n-1$ the matrix $B_1^kA_2$ is very proximal.

Let $\overrightarrow{u}_1\in\tilde{{\mathbb{K}}}_v\cdot\hat{u}_1$
be an eigenvector of $A_1$ corresponding to $\alpha_1$ as in
Lemma \ref{nice-vectors}, i.e. the coordinates of $\overrightarrow{u}_1$ are $%
\tilde{S}$-integers, and $|\overrightarrow{u}_1|_w\leq
|\alpha_1/\alpha_2|_v^{r_4}$, for any $w\in\ti S$.

For any $k\leq 2n-1$ we have
\begin{equation*}
\| B_1^k(\overrightarrow{u}_1)\|_w\leq \| B_1\|_w^{2n-1}\|\overrightarrow{u}%
_1\|_w\leq |\alpha_1/\alpha_2|_v^{(2n-1)m_2r_0+r_4},
\end{equation*}
for any $w\in\ti S$, while for any $w\notin\ti S$ the $w$-norm of this vector is $%
\leq 1$. It follows that for any $1\leq k_1<\ldots<k_n\leq 2n-1$
we have
\begin{equation*}
|\text{det}(B_1^{k_1}(\overrightarrow{u}_1),\ldots,B_1^{k_n}(\overrightarrow{u}%
_1))|_w \leq |\alpha_1/\alpha_2|_v^{((2n-1)m_2r_0+r_4)n}
\end{equation*}
for any $w\in\ti S$, and
\begin{equation*}
|\text{det}(B^{k_1}(\overrightarrow{u}_1),\ldots,B^{k_n}(\overrightarrow{u}%
_1)|_w \leq 1
\end{equation*}
for any $w\notin\ti S$. Since $B_1$ is in $(2n-1)$--general
position with respect to the $\{\overrightarrow{u}_i\}$'s, this
determinant is not zero, and hence by the product formula
\begin{equation*}
\prod_{\text{{\tiny over all places}}} |\text{det}(B_1^{k_1}(%
\overrightarrow{u}_1),\ldots,B_1^{k_n}(\overrightarrow{u}_1)|_w
=1,
\end{equation*}
which implies that
\begin{equation*}
\prod_{w\in\ti S} |\text{det}(B_1^{k_1}(\overrightarrow{u}_1),\ldots,B_1^{k_n}(%
\overrightarrow{u}_1))|_w \geq 1.
\end{equation*}
We conclude:

\begin{cor}
For any $w\in\ti S$, and in particular for $w=v$
\begin{equation*}
|\text{det}(B_1^{k_{1}}(\overrightarrow{u}_{1}),\ldots ,B_1^{k_{n}}(%
\overrightarrow{u}_{1}))|_{w}\geq |\alpha _{1}/\alpha
_{2}|_{v}^{-((2n-1)m_2r_{0}+r_{4})n|\ti S|}.
\end{equation*}
\end{cor}

We will need also the following:

\begin{lem}
\label{v-f-H} Suppose that $\overrightarrow{v}_{1},\ldots ,\overrightarrow{v}%
_{n}$ are any $n$ vectors in $\tilde{{\mathbb{K}}}_{v}^{n}$
satisfying

\begin{itemize}
\item  $\Vert \overrightarrow{v}_{i}\Vert _{v}\leq t^{c^{\prime
}},~\forall i\leq n$, and

\item  $|\text{det}(\overrightarrow{v}_{1},\ldots ,\overrightarrow{v}%
_{n})|_{v}\geq t^{-c^{\prime \prime }}$,
\end{itemize}

for some $c^{\prime },c^{\prime \prime }\in {\mathbb{N}}$ and
$t>0$.

Then for any hyperplane $H\subset \tilde{\BK}^{n}_v$ there is
$i\leq n
$ such that $d([\overrightarrow{v}_{i}],[H])\geq \frac{1}{\lambda_{1}\lambda _{n-1}}%
t^{-c^{\prime \prime }-(n-1)c^{\prime }}$ in the $\tilde{\BK}_v$
projective space, where $\lambda_{k}$ is the volume of the
$k$--dimensional unit ball (in particular $\lambda _{k}=1$ in the
non-Archimedean case).
\end{lem}

\begin{proof}
Let $f$ be a linear form such that $\|f\|=1$ and $H=\ker(f)$, then
the volume of $\{x\in\tilde{\BK}_v^n : |f(x)|_v \leq |a|_v,
\|x\|_v \leq |b|_v\}$ is bounded above by
$\lambda_{1}\lambda_{n-1}|a|_v|b|_v^{n-1}$ -- the volume of a
``cylinder'' with base radius $|b|_v$ and ``height'' $2|a|_v$, for
any $a,b \in \tilde{\BK}_v$. Since
$d([x],[H])=\frac{|f(x)|_v}{\|x\|_v}$, we get the desired
conclusion by comparing this volume to the determinant of the
${\overrightarrow{v}_{i}}$'s.
\end{proof}

Setting $c^{\prime}=(2n-1)m_2r_0+r_4$ and
$c^{\prime\prime}=((2n-1)m_2r_0+r_4)n|\ti S|$ we
get some constant\footnote{%
Note that we can fix $r_{8}$ before determining $r_6,r_7$.} $r_8$,
such that whenever
$\overrightarrow{v}_1,\ldots,\overrightarrow{v}_n$ are as in Lemma
\ref{v-f-H} with $t=|\alpha_1/\alpha_2|_v$ and $[H]$ is some
projective
hyperplane, there is one $[\overrightarrow{v}_i]$ at distance at least $%
|\alpha_1/\alpha_2|^{-r_8}_v$ from $[H]$, in particular:

\begin{lem}
\label{r10}\label{const:r8} For any $1\leq k_{1}<k_{2}<\ldots
<k_{n}\leq 2n-1 $ and any hyperplane $H\subset \tilde{\BK}^{n}_v$
there exists $i\leq n $ such that
\begin{equation*}
d([B_1^{\pm k_{i}}\hat{u}_{1}],[H])\geq |\alpha _{1}/\alpha
_{2}|_{v}^{-r_{8}}.
\end{equation*}
\end{lem}

By the pigeonhole principle, we conclude:

\begin{cor}\label{pigeon}
\label{2n-1} For any two hyperplanes $H_{1},H_{2}\subset
\tilde{\BK}^{n}_v$ there is some $k\leq 2n-1$ such that we have
simultaneously
\begin{equation*}
d([B_1^{k}\hat{u}_{1}],[H_{1}])\geq |\alpha _{1}/\alpha
_{2}|_{v}^{-r_{8}}, d([B_1^{-k}\hat{u}_{1}],[H_{2}])\geq |\alpha
_{1}/\alpha _{2}|_{v}^{-r_{8}}.
\end{equation*}
\end{cor}

Now let $[H^+],[H^-]$ be the repelling hyperplanes for the $%
|\alpha_1/\alpha_2|_v^{-r_6}$--very contracting element $%
A_2=A_1^{r_7}B_1A_1^{-r_7}$ and its inverse, and take the
corresponding $k$ in Corollary \ref {2n-1}. Recall that the
attracting points $t^+,t^-$ of $A_2^{\pm 1}$ are both at distance
at most $|\alpha_1/\alpha_2|_v^{-r_6}$ from $[\hat u_1]$. We thus
obtain:

\begin{prop}[The very proximal element $X$]\label{prop:very-proximal}
Assume that $r_{8}>2(2n-1)r_{0}$, then the element
$X=B_1^{k}A_{2}$ is $(\rho,\gd)$--very proximal with
$$
 \rho=|\alpha_{1}/\alpha _{2}|_{v}^{-2m_2kr_{0}}(|\alpha _{1}/
 \alpha_{2}|_{v}^{-r_{8}}-|\alpha _{1}/
 \alpha_{2}|_{v}^{-r_{6}+4m_2kr_{0}}),\text{~and~}
 \gd=|\alpha_{1}/\alpha _{2}|_{v}^{-r_{6}+4m_2kr_{0}}
$$
and with repelling hyperplanes
\begin{equation*}
\lbrack H_{X}^{+}]=[H^{+}],~~[H_{X}^{-}]=B_1^{k}[{H}^{-}]
\end{equation*}
and attracting points
\begin{equation*}
\lbrack t_{X}^{+}]=B_1^{k}t^{+},~~[t_{X}^{-}]=t^{-}.
\end{equation*}
\end{prop}

\begin{proof}
Since $\|B_1^{\pm 1}\|_v\leq |\ga_1/\ga_2|_v^{m_2r_0}$, $B_1$ is
$|\ga_1/\ga_2|_v^{4m_2r_0}$ bi-Lipschitz on the entire projective
space. This implies that $X=B_1^kA_2$ is
$|\ga_1/\ga_2|_v^{-r_6+4m_2kr_0}$ very contracting with the
specified attracting points and repelling hyperplanes, and that
\begin{eqnarray*}
 d(B_1^k (t^+),[H^+])\geq
 d(B_1^k[\hat u_1],[H^+])-d(B_1^k [\hat
 u_1],B_1^k (t^+))\geq
 |\ga_1/\ga_2|_v^{-r_8}-|\ga_1/\ga_2|_v^{-r_6+4km_2r_0},
\end{eqnarray*}
and
\begin{eqnarray*}
 d(t^{-},B_1^{k}[H^{-}])\geq\|B_1^{\pm k}\|_v^{-4}d(B_1^{-k}
 (t^{-}),[H^{-}])
 \geq
 |\ga_1/\ga_2|_v^{-4m_2kr_0}(|\ga_1/\ga_2|_v^{-r_8}-|\ga_1/\ga_2|_v^{-r_6+4km_2r_0})
\end{eqnarray*}
\end{proof}

Taking $r_6>>r_8$ sufficiently large (after choosing $r_8$
sufficiently large) we may assume that:
\begin{equation*}
\rho=|\alpha_1/\alpha_2|_v^{-r_8-2m_2kr_0}(1-|\alpha_1/\alpha_2|_v^{r_8-r_6+6m_2kr_0})%
\geq\frac{1}{2}|\alpha_1/\alpha_2|_v^{-2r_8}.
\end{equation*}
Set
\begin{equation*}
r_9=2r_8+d\log_\gt 2,~~r_{10}=r_6-4m_2kr_0,
\end{equation*}
Then we get that $X=B_1^kA_1^{r_7}B_1A_1^{-r_7}$ is $(|\alpha_1/%
\alpha_2|_v^{-r_9},|\alpha_1/\alpha_2|_v^{-r_{10}})$--very
proximal. The matrix $X$ is our first ping--pong player.

%-----------------------------------------------------------------------------------------------------------------------------------

\subsection{Step 4}\label{sbs:step4}

The last step of the proof consists in finding a second ping--pong
partner $Y$ by conjugating $X$ by a suitable bounded word in the alphabet $\gS$. This is performed in quite the same way as in Step 3, so
we only sketch the proof here.

Note first that $X$ is a word in $\gS$ of length at most
$2(2n-1)m_2+2r_7r_3r$. Therefore, by
requiring $s$ from Proposition \ref{sigma-two} to be at least this constant, we can
assume that $X$ is semisimple. Let $[\hat v_1]$ (resp. $[\hat
v_n]$) be the eigendirection of the maximal (resp. minimal)
eigenvalue of $X$.

Let $B_2$ be a word in $\gS$ which is in $(2n-1)^2$--general
position with respect to the eigenvectors of $X$ (chosen as in
Lemma \ref{nice-vectors}). Again by Lemma \ref{Bezout} and the
discussion following Definition \ref{defn:general-position}, we
may find $B_2$ as a word of length $\leq m_2((2n-1)^2)$.
We can then apply the same pigeonhole argument as in Corollary
\ref{pigeon} and obtain an index $k'\leq (2n-1)^2$ such that
$B_2^{k'}[\hat v_1]$ and $B_2^{k'}[\hat v_n]$ are both far away
from the repelling hyperplanes $[H_X^+]$ of $X$ and $[H_X^-]$ of
$X^{-1}$ (i.e. $|\alpha_1/\alpha_2|_v^{-r_{11}}$--apart for some
other constant $r_{11}$).

Setting $Y=B_2^{k'}XB_2^{-k'}$, we see that some bounded power
$Y^{r_{12}}$ of $Y$ is very proximal with attracting and repelling
points $B_2^{k'}[\hat v_1]$ and $B_2^{k'}[\hat v_n]$ and repelling
hyperplanes $B_2^{k'}[H_X^+]$ and $B_2^{k'}[H_X^-]$. Since those
points are away from the repelling hyperplanes of $X$ (or any
power of $X$), we conclude that, after taking a larger power
$r_{13}\geq r_{12}$ if necessary, $Y^t$ and $X^{t}$
play ping--pong, and hence independent, for any $t\geq r_{13}$.

\begin{rem}
As mentioned at the beginning of Section \ref{sec:ping-pong} we
assumed throughout that the representation space $V$ was
irreducible. Lemma \ref{nice-vectors} as well as the rest of the
argument above, relies on the assumption that the entries of the
elements of $\gS$ viewed as matrices acting on $V$ are
$S$--arithmetic. However, in general our wedge representation $V$
might be reducible, and we have to replace it with some
irreducible subquotient where this assumption may not hold. In
order to cope with this problem we argue as follows. We change the
representation space from $V$ to an irreducible subquotient
$V_0/W$ where $V_0,W$ are invariant subspaces of $V$. One can
carry out the proof of Lemma \ref{nice-vectors} in $V$ and first
treat the eigenvectors in $W$, then those in $V_0\setminus W$ and
finally take the projections of those to $V_0/W$. This would yield
an analogous statement for $V_0/W$ which is sufficient for the
whole argument. Note also that in characteristic zero, as
$\BG$ is semisimple, $V$ is completely reducible so our
irreducible representation is a sub-representation of the wedge
power, rather than a subquotient, and hence, in this case, we may take $V_0$
instead of $V$ without further changes.
\end{rem}

This completes the proof of Theorem \ref{thm2}. \qed

\medskip

\section{A Zariski dense free subgroup in characteristic zero}\label{sec:Z-dense}

We will now give two stronger versions of Theorem
\ref{thm1} which are useful for applications. Since all the applications we have in mind are for fields of characteristic zero, we allow ourselves to make this restriction, although we believe that it is unnecessary.

\begin{thm}\label{thm:Z-dense}
Let $K$ be a field of characteristic zero, $\BH$ a Zariski
connected semisimple $K$--group and $\gC\leq H=\BH(K)$ a finitely
generated Zariski dense subgroup. Then there is a constant
$m_1=m_1(\gC)$ such that for any symmetric generating set $\gS\ni
1$ of $\gC$, $\gS^{m_1}$ contains two independent elements which
generate a Zariski dense subgroup of $\BH$.
\end{thm}

\begin{rem}\label{rem:Z-dense}
The proof of Theorem \ref{thm:Z-dense} also shows that Theorem \ref{thm2} remains true, in characteristic zero, with the stronger conclusion that the independent
elements generate a Zariski dense subgroup of $\BG$.
\end{rem}

In order to obtain Theorem \ref{thm:Z-dense} one needs to slightly
modify the argument of Section \ref{sec:ping-pong} in a few places.
We will now indicate these modifications. For the sake of
simplicity, let us assume that $\BH$ is simple.

It is well known that $\BH$ admits two conjugate elements which
generate a Zariski dense subgroup. Indeed, one can take a regular
unipotent in $\BH$ and a conjugate lying in an opposite
parabolic (these unipotent elements can be taken in $\BH(\tilde{K})$,
where $\tilde{K}$ is a finite extension of $K$ over which $\BH$ is
isotropic). Let $\mathcal{A}$ be the subalgebra spanned by
$\text{Ad}(\BH)$ in $\text{End}(\mathfrak{h})$ where
$\mathfrak{h}$ denotes the Lie algebra of $\BH$, set
$$
 F=\{ (g,h)\in\BH\times\BH:\text{Ad}(g)~\text{and}~\text{Ad}(h)~\text{do not generate the algebra}~\mathcal{A}\},
$$
and $E=\{ (g,h)\in\BH\times\BH:(g,hgh^{-1})\in F\}$.
It follows the algebraic variety $E$
is proper in $\BH\times\BH$. Let $E_1=\{ g\in \BH: (g,h)\in
E~\forall h\in\BH\}$, and for $g\in\BH$ let $E_2(g)=\{
h\in\BH:(g,h)\in E\}$. Then $E_1$ is a proper subvariety of $\BH$
and one easily checks that $\chi (E_2(g))$ is bounded independently
of $g$.

Note that in Theorem \ref{thm:Z-dense}, we did not assume that the field is a global field. In order to take
 care of this issue we will use the specialization map introduced in Lemma \ref{lem:spec}. Without loss of generality,
 we may pass to a subgroup of finite index in $\gC$. We will need the following lemma.

\begin{lem}\label{graphdensity}
Let $f: \Gamma \mapsto \BG$ be the specialization map from Lemma \ref{lem:spec}. Then the subgroup $\gD=\{(\gc, f(\gc)) \in \BH \times \BG\ / \gc \in \gC\}$ is not contained in any algebraic subset of the form $(V \times \BG) \cup (\BH \times W)$, where $V$ and $W$ are proper closed subvarieties of $\BH$ and $\BG$ respectively.
\end{lem}

\begin{proof}
This is obvious since $\gC$ is Zariski dense in $\BH$ and $f(\gC)$
is Zariski dense in $\BG$.
\end{proof}

When pursueing the argument of Section 6, we need to specify conditions on the elements of the generating set. These conditions are set on elements of $f(\gC) \in \BG$. We will now introduce new algebraic conditions directly on the elements of $\gC \in \BH$. Combining Lemma \ref{Bezout} with Lemma \ref{graphdensity} we see that given a set of non-trivial algebraic conditions in $\BH$ and another such set in $\BG$, there is an integer $N$ such that, for every generating set $\gS$ of $\gC$, there is a point $\gc \in \gS^N$ such that $\gc\in \BH$ and $f(\gc) \in \BG$ do not satisfy those conditions. This will be used repeatedly below. When no confusion may arise we will often say for instance that some element $A \in \gC$ acts proximally when we really mean that $f(A) \in \BG$ acts proximally on the representation variety used in Section 6.

The first modification needed in the argument of Section
\ref{sec:ping-pong} is in Proposition \ref{prop:very-proximal}
when we construct the very proximal element $X$. Instead of using
the same element $B_1$ which was used in the construction of the
very contracting element, we should use an element $B_1'$ which
satisfies
\begin{itemize}
\item $f(B_1')$ is in $(2n-1)$--general position with respect to $\{\overrightarrow{u_i}\}_{i-1}^n$ (like $f(B_1)$), and
\item $(B_1')^k\notin E_1A_2^{-1},~\forall k\leq 2n-1$.
\end{itemize}
Note that the choice of $B_1'$ depends on $A_2$, however, since
$\chi (E_1A_2^{-1})$ is independent of $A_2$ we can find $B_1'$ in a
fixed power of our generating set $\gS$. Retrospectively we should
also take the constant $r_6$ big enough so that the very
contracting element $A_2$ constructed in Proposition
\ref{prop:very-proximal} will have sufficiently small attracting
and repelling neighborhoods (i.e. that $|\ga_1/\ga_2|^{-r_6}$ will
be small enough) so that the element $X=(B_1')^kA_2$ (where $k$ is
some integer $\leq 2n-1$) becomes very proximal. Additionally, we
have to take the constant $s$ in Proposition \ref{sigma-two}
sufficiently large to guarantee that $f(X)$ is still semisimple.

The second change one has to do is in Step (4) when choosing the
appropriate conjugation of $X$. By the choice of $B_1'$ we know
that $X\notin E_1$. We take $B_2'$ which satsfies:
\begin{itemize}
\item $f(B_2')$ is in $(2n-1)^2$--general position with respect to the eigenvectors of $f(X)$ (again chosen as in Lemma \ref{nice-vectors}), and
\item $(B_2')^k\notin E_2(X),~\forall k\leq (2n-1)^2$.
\end{itemize}
Then, as in the previous section, if $Y=(B_2')^{k'}X(B_2')^{-k'}$
for some appropriate $k'\leq (2n-1)^2$ then $X^t$ and $Y^{t}$ are
independent for any $t\geq r_{13}'$ for some constant $r_{13}'$.

Finally, since $F$ is an algebraic subvariety of $\BH\times\BH$ and $\chi(Fx)$ is bounded independently of $x \in \BH\times\BH$, we may apply
Lemma \ref{Bezout} to the set $\{(X,Y)\}$ and the variety $F\cdot (X^{-r_{13}'}, Y^{-r_{13}'})$. Since $\{(X,Y)\}$ is not in $F$, it follows that $\{(X,Y)\}$ generates a group not contained in $F\cdot (X^{-r_{13}'}, Y^{-r_{13}'})$. Hence Lemma \ref{Bezout} yields some $t$ with $r_{13}'\leq t\leq r_{13}'+N(\chi(F))$ such that $(X^t,Y^t)\notin F$. Now since $X$ has infinite order, the
Zariski connected group $(\overline{\langle
X^t,Y^t\rangle}^Z)^\circ$ has positive dimension, and since it is
normalized by $X^t$ and $Y^t$, while
$\text{span}\{\text{Ad}(X^t),\text{Ad}(Y^t)\}=\mathcal{A}$ it
follows that $(\overline{\langle X^t,Y^t\rangle}^Z)^\circ$ is
normal in $\BH$. Since $\BH$ is assumed to be simple we derive
that $\langle X^t,Y^t\rangle$ is Zariski dense.\qed

\begin{thm}\label{thm:strong-version}
Let $\BG$ be a semisimple algebraic group defined over a field $K$ of characteristic zero, $\gC$ a finitely generated Zariski dense subgroup of
$\BG(K)$, and $V\subset\BG\times\BG$ a proper algebraic
subvariety. Then there is a constant $m=m(\gC,V)$ such that for
any generating set $\gS\ni 1$ of $\gC$, $\gS^m$ contains a pair of
independent elements $x,y$ with $(x,y)\notin V$.
\end{thm}

\begin{proof} The subset $\gS^{m_1}$ contains a pair $\{A,B\}$ of independent elements for some constant $m_1=m_1(\gC)$ given by Theorem \ref{thm:Z-dense}. This pair generates a Zariski dense subgroup of $\gC$. Hence $(1,A)$, $(1,B)$, $(A,1)$ and $(B,1)$ together generate a Zariski dense subgroup of $\BG \times \BG$. The set $V'=V \cup \{(x,y) | [x,y]=1\}$ is a proper closed algebraic subset of $\BG \times \BG$. By Lemma \ref{Bezout}, there exists another constant $m_2=m_2(V)$ such that
some word of length at most $m_2$ in those four generators lies
outside $V'$. This word has the form $(W_1(A,B),W_2(A,B))$ where
$W_1,W_2$ are bounded words in $A$ and $B$ that do not commute as
words in the free group. It follows that they generate a free
subgroup, hence form a pair of independent elements in
$\gS^{m_1m_2}$.
\end{proof}

%\begin{remark}
%In the above theorem, we can actually assume that $\gC$ is a
%subgroup of  $\BG(K) \leq \GL_d(K)$, where $K$ is an arbitrary
%field, as opposed to only a global field. Indeed one can always
%map $\gC$ to $\GL_d(\OO_\BK (S))$, where $\OO_\BK (S)$ is the ring
%of $S$--integers of some global field $\BK$, in such a way that
%the image is not virtually solvable (this is Lemma
%\ref{lem:spec}). Moreover, it can be shown, by the same argument
%of specialization used in \cite{EMO}, that the ring $\OO_\BK (S)$
%depends only on the ring generated by the matrix entries of the
%elements of $\gC$. We can thus apply Proposition \ref{sigma-two}
%to $\gC$ in $\BG$ with the number $s$ being taken to be
%$s=m(d,\BK,S)$, the constant from Theorem \ref{thm2}. We get a
%pair $\{A,B\}$ and can then specialize the group $\langle A,B
%\rangle$ generated by this pair and apply Theorem \ref{thm2} to
%its image under the specialization.
%\end{remark}

%----------------------------------------------------------------------------------

\section{Some applications}\label{sec:applications}

In this section we draw some consequences of our main result.

\subsection{Uniform non-amenability and a uniform Cheeger
constant.}\label{subsec:u-n-a}

Recall that a group is called amenable if the regular
representation admits almost invariant vectors. It follows that if
a non-amenable group $\gC$ is generated by a finite set $\gS$ then
there is a positive constant $\gep(\gS )$ such that for any $f\in
L^2(\gC )$ there is some $\gs\in\gS$ for which $\|\rho (\gs
)(f)-f\|\geq\gep (\gS )\| f\|$, where $\rho$ denotes the left
regular representation, i.e. $\rho(\gc )(f)(x):=f(\gc^{-1}x)$.
Such an $\gep (\gS )$ is called a {\it Kazhdan constant for
$(\gS,\rho )$}. A finitely generated group $\gC$ is said to be
{\it uniformly non-amenable} if there is a positive Kazhdan
constant $\gep =\gep(\gC )>0$ for the regular representation which
is independent of the generating set $\gS$, i.e. if there is
$\gep>0$ such that for any generating set $\gS$ of $\gC$ and any
$f\in L^2(\gC )$ there is $\gs\in\gS$ for which $\|\rho (\gs
)(f)-f\|\geq \gep\| f\|$. It was observed by Y. Shalom \cite{Sha}
that Theorem \ref{thm1} implies:

\begin{thm}\label{thm:u-n-a}
A finitely generated non-amenable linear group is uniformly
non-amenable.
\end{thm}

\begin{proof}
The proof is an elaboration of the original proof by Von-Neumann
that a group which contains a non-abelian free subgroup is
non-amenable.

Let $\gC$ be a non-amenable linear group, and let $m=m(\gC )$ be
the constant from Theorem \ref{thm1}. Let $\gS$ be a generating set
of $\gC$ and let $x,y\in (\gS\cup\gS^{-1}\cup 1)^m$ be two
independent elements. Denote by $F_2=\langle x,y\rangle$ the
corresponding free subgroup. Choose a complete set $\{ c_i\}$ of
right coset representatives for $F_2$ in $\gC$, and write
$$
 L^2(\gC )=\bigoplus L^2(F_2c_i).
$$
Let $f\in L^2(\gC )$, and let $f_i$ denote the restriction of $f$
to $F_2c_i$. Let $\gt_0$ be the Kazhdan constant for
$(\rho_{F_2},\{x,y\})$ then for any $i$ either $x$ or $y$ moves
$f_i$ by at least $\gt_0\| f_i\|$. Let
$$
 f_x=\sum_{\|\rho (x)f_i-f_i\|\geq\gt_0\| f_i\|}f_i, ~~\text{and}~~
 f_y=\sum_{\|\rho (y)f_i-f_i\|\geq\gt_0\| f_i\|}f_i
$$
Then either $\| f_x\|\geq\|f\|/\sqrt{2}$ or $\|
f_y\|\geq\|f\|/\sqrt{2}$. Without loss of generality let us assume
that $\| f_x\|\geq{\| f\|}/{\sqrt 2}$. It follows that
$$
 \|\rho (x)f-f\|\geq\|\rho (x)f_x-f_x\|\geq\gt_0\|
 f_x\|\geq\frac{\gt_0}{\sqrt{2}}\|f\|.
$$
Now write $x=\gs_1^{\gep_1}\cdot\ldots\cdot\gs_m^{\gep_m}$ where
$\gs_i\in\gS\cup\{ 1\}$ and $\gep_i=\pm 1$. Then by the triangle
inequality, if we let $\gs_0=1$ and $\gep_0=1$
$$
 \|\rho (x)f-f\|\leq\sum_{i=1}^{m}\|\rho
 (\gs_0^{\gep_0}\cdot\ldots\cdot\gs_i^{\gep_i})f-\rho
 (\gs_0^{\gep_0}\cdot\ldots\cdot\gs_{i-1}^{\gep_{i-1}})f\|=
 \sum_{i=1}^m\|\rho (\gs_i)f-f\|,
$$
and hence, for some $i$ we have $\|\rho(\gs_i
)f-f\|\geq\frac{\gt_0}{m\sqrt{2}}\| f\|$.
\end{proof}

Note that usually such groups do not admit a uniform Kazhdan
constant for arbitrary unitary representation, even if they have
property $(T)$ (see \cite{GZ}).

By considering $f$ to be a characteristic function of a finite
subset of $\gC$ and applying Theorem \ref{thm:u-n-a} we obtain the
following useful result. We denote by $|A|$ the number of elements
in $A$ and by $\bigtriangleup$ the operator of symmetric
difference between sets.

\begin{thm}\label{cor:boundary}
Let $\gC$ be a finitely generated non-virtually solvable linear
group. Then there is a positive constant $b=b(\gC )$ such that for
any generating set $\gS$ (not necessarily finite or symmetric) of
$\gC$, and any finite subset $A\subset\gC$ there is $\gs\in\gS$
such that:
$$
 \frac{|\gs\cdot A\bigtriangleup A|}{|A|}\geq b.
$$
\end{thm}

Consider a graph $X$ and a finite subset $A\subset X$. The {\it
boundary} of $A$ is the set $\partial A$ of all vertices in $A$
which have at least one neighbor outside $A$. The {\it Cheeger
constant} $\mathcal{C}(X)$ is defined by
$$
 \mathcal{C}(X)=\inf\frac{|\partial A|}{|A|},
$$
where $A$ runs over all finite subsets of $\text{vert}(X)$ when
$X$ is infinite, and over all subsets of size at most
$|\text{vert}(X)|/2$ when $X$ is finite. For a group $\gC$ and a
finite generating set $\gS$ we denote by $\mathcal{C}(\gC,\gS)$
the Cheeger constant of the Cayley graph of $\gC$ with respect to
$\gS$, and by $\mathcal{C}(\gC )$ the {\it uniform Cheeger
constant} of $\gC$:
$$
 \mathcal{C}(\gC ):=\inf\{\mathcal{C}(\gC,\gS):\gS\text{~is a finite generating
 set~}\}.
$$
In some places (c.f. \cite{Arz},\cite{Sha}) a group is called
uniformly non-amenable if it has a positive uniform Cheeger
constant. Clearly our definition of uniform non-amenability
implies this one\footnote{It is still not known whether these two
definitions are equivalent for a general finitely generated group.}:

\begin{cor}\label{cor:expanders}
A finitely generated non-virtually solvable linear group has a
positive uniform Cheeger constant.
\end{cor}

When $\mathcal{C}(X)>\gep$ the graph $X$ is said to be {\it
$\gep$--expander}. Hence Corollary \ref{cor:expanders} can be
reformulated as follows:

\begin{cor}
The family of all Cayley graphs corresponding to finite generating
sets of a given non-virtually solvable linear group $\gC$ form a
family of $\gd$--expanders for some constant $\gd=\gd (\gC )>0$.
\end{cor}

%---------------------------------------------------------------------------

\subsection{Growth}\label{sub-sec:growth}

In \cite{EMO}, Eskin Mozes and Oh proved that any finitely
generated non-virtually solvable linear group has a uniform
exponential growth by showing that some bounded words in the
generators generate a free semigroup. As a consequence of Theorem
\ref{thm1} (more precisely of Theorem \ref{cor:boundary}) we
obtain:

\begin{thm}\label{thm:growth}
Let $\gC$ be a finitely generated non-virtually solvable linear
group. Then there is a constant $\lambda =\gl (\gC )>1$ such that
if $\gS$ is
%any set which generate a group whose Zariski closure
%contains the identity connected component of the Zariski closure
%of $\gC$,
any finite generating set of $\gC$, then $|\gS^n|\geq |\gS |\gl^{n-1},~\forall n\in\BN$.
\end{thm}

Since the proof is straightforward, we will omit it. One can
actually take $\gl =1+\frac{b}{2}$ where $b$ is the constant from
Theorem \ref{cor:boundary}.

\begin{rem}
Theorem \ref{thm:growth} improves Eskin-Mozes-Oh theorem in
several aspects:
\begin{itemize}
\item Unlike the situation in \cite{EMO}, the generating set $\gS$
in Theorem \ref{thm:growth} is not assumed to be symmetric, so
\ref{thm:growth} gives the uniform exponential growth for
semigroups rather than just for groups.

\item We didn't make the assumption from \cite{EMO} that the
characteristic of the field is $0$.

\item The estimate on the growth that we obtain is sharper; In
particular if the generating set is bigger the growth is faster.
This sharper estimate is important for applications,
% e.g. it implies that manifolds with almost non-negative Ricci curvature
%has almost nilpotent fundamental group (in the case that the
%fundamental group is linear) \cite{Petrunin}.
\end{itemize}
\end{rem}

As another consequence we obtain that also the spheres have
uniform exponential growth, and moreover the size of each sphere
is at least $\frac{b}{2}$ times the size of the corresponding
ball. If $\gS\ni 1$ is a generating set for $\gC$ the sphere
$S(n,\gS )$ corresponding to $\gS$ is the set of all elements in
$\gC$ of distance exactly $n$ from $1$ in the Cayley graph, i.e.
$S(n,\gS )=\gS^n\setminus\gS^{n-1}$.

\begin{cor}
Let $\gC,\gS$ and $\gl=1+\frac{b}{2}$ be as in Theorem
\ref{thm:growth}. Then
$$
 |S(n,\gS)|\geq\frac{b}{2}|\gS^{n-1}|\geq
 \frac{b}{2}|\gS |(1+\frac{b}{2})^{n-2}.
$$
\end{cor}

One can derive many other variants of these results from Theorem
\ref{cor:boundary}. Here is another example:

\begin{exercise}
Let $\gC,\gS$ and $\gl$ be as above. There is a sequence
$\{\gs_i\}_{i \in \BN}$ of elements of $\gS$ such that for any
$n\in\BN$
$$
 |\{\prod_{1\leq i_1<\ldots<i_k\leq
 n}\gs_{i_1}\cdot\ldots\cdot\gs_{i_k}\}|\geq\gl^{n-1}.
$$
\end{exercise}

%------------------------------------------------------------------------------

\subsection{Dense free subgroups, amenable actions and growth of
leaves.}\label{seb-sec:dense}

Theorem \ref{thm1} implies the following result from \cite{BG}
which answered a question of Carri\`{e}re and Ghys \cite{Ghys}:

\begin{thm}\label{thm:dense}
Let $G$ be a connected semisimple Lie group and $\gC\leq G$ a
dense subgroup. Then $\gC$ contains a dense free subgroup of rank
$2$.
\end{thm}

\begin{proof}
Let us assume for simplicity that $G$ is simple. The adjoint
representation $\text{Ad}:G\to\GL(\mathfrak{g})$ is irreducible
and by Burnside's theorem its image spans
$\text{End}(\mathfrak{g})$. It is well known that
$\text{End}(\mathfrak{g})$ is generated by two elements and that
these elements can be chosen in $\text{Ad}(G)$. Since
$\text{End}(\mathfrak{g})$ is finite dimensional it follows that
the set
$$
 V=\{ (g_1,g_2)\in G\times G:
 \text{Ad}(g_1)\text{~and~}\text{Ad}(g_2)\text{~generate~}\text{End}(\mathfrak{g})\}
$$
is Zariski open $G\times G$. By Theorem \ref{thm:strong-version}
and the remark following it there is a constant $m=m(\gC,V)$ such that if $\gS\ni 1$ is any
generating set for $\gC$ then $\gS^m$ contains independent
elements $x,y$ with $(x,y)\notin V$. Let $U\subset G$ be a
Zassenhaus neighborhood (c.f. \cite{raghunathan} 8.16), and let
$\gO$ be an identity neighborhood with $\gO^m\subset U$. Take
$\gS=\gC\cap\gO$. Since $G$ is connected and $\gC$ is dense, $\gS$
generates $\gC$, and therefore $\gS^m$ contains $x,y$ independent
with $(x,y)\notin V$ according to Theorem
\ref{thm:strong-version}. Now the connected component of the
identity in $\overline{\langle x,y\rangle}$ is normalized by $x,y$
and as $(x,y)\notin V$ it is normal in $G$, and by simplicity of
$G$ it is either $1$ or $G$. In other words $\langle x,y\rangle$
is either discrete or dense. However $\langle x,y\rangle$ is free
and hence non-nilpotent and since $x,y\in U$ it follows from
Zassenhaus' theorem that $\langle x,y\rangle$ is not discrete.
\end{proof}

Recall that one of the main motivation to prove Theorem
\ref{thm:dense} was the Connes--Sullivan conjecture which was
first proved by Zimmer:

\begin{cor}[Zimmer \cite{Zim}]
Let $\gC$ be a countable subgroup of a real Lie group $G$. Then
the action of $\gC$ on $G$ by left translations is amenable iff
the connected component of the identity of the closure of $\gC$ is
solvable.
\end{cor}

This corollary follows from Theorem \ref{thm:dense} by the
observation of Carri\`{e}re and Ghys that a non-discrete free
subgroup of $G$ cannot act amenably (see \cite{BG},\cite{BG1} for
more details and stronger results).

Another motivation was the result about the
polynomial--exponential dichotomy for the growth of leaves in
Riemannian foliations which was conjectured by Carri\`{e}re:

\begin{thm}[\cite{BG1}]\label{thm:foliations}
Let $\mathcal{F}$ be a Riemannian foliation on a compact manifold.
Then either the growth of any leaf in $\mathcal{F}$ is polynomial,
or the growth of a generic leaf is exponential.
\end{thm}

Theorem \ref{thm:foliations} can be considered as a foliated
version of the well known conjecture according to which the growth
of the universal cover of any compact Riemannian manifold is
either polynomial or exponential. The proof of
\ref{thm:foliations} relies on the following strengthening of
Theorem \ref{thm:dense} as well as some special argument for
solvable groups (see \cite{BG1} for more details):

\begin{thm}(\cite{BG})\label{thm:dense1}
Let $G$ be a connected semisimple Lie group, $\gC\leq G$ a dense
subgroup and $\gO_1,\ldots,\gO_n$ some $n$ open sets in $G$. Then
one can pick $x_i\in\gC\cap\gO_i,~i=1,\ldots,n$ which are
independent, i.e. generate a free group of rank $n$.
\end{thm}

In \cite{BG} Theorem \ref{thm:dense1} was the main result and
Theorem \ref{thm:dense} followed as a consequence. Let us show
that conversely it is possible to derive Theorem \ref{thm:dense1}
from Theorem \ref{thm:dense}. This way, Theorem \ref{thm:dense1}
will appear as a mere consequence of the main theorem of the
present paper, namely Theorem \ref{thm1}.

\textit{Proof that Theorem \ref{thm:dense} implies Theorem
\ref{thm:dense1}}. Let $F_2\leq\gC$ be a free subgroup of rank $2$
which is dense in $G$, and let $F_n$ be a subgroup of index $n-1$
in $F_2$. Then $F_n$ is a free group of rank $n$ which is still
dense in $G$. We will pick the $x_i$ in $F_n$ inductively as
follows. Suppose we picked already $x_1,\ldots,x_{i-1}$. Since
$F_n$ is dense and $\gO_i$ is open, $F_n$ is generated by
$F_n\cap\gO_i$. It follows that we can pick $x_i\in F_n\cap\gO_i$
such that the abelianization of $\langle x_1,\ldots,x_i\rangle$
has rank $i$; Indeed, look at the tensor of the abelanization of
$F_n$ with $\BQ$ and pick $x_i$ in the generating set
$F_n\cap\gO_i$ which is not in the $(i-1)$--dimensional
$\BQ$--subspace spanned by the images of $x_1,\ldots,x_{i-1}$. It
follows that $\langle x_1,\ldots,x_i\rangle$ is a free group whose
minimal number of generators is exactly $i$. Since a free group is
Hopfian it follows that $x_1,\ldots,x_i$ are independent. \qed

\medskip

We refer the reader to \cite{BG1} for an extension of
Theorems \ref{thm:dense} and \ref{thm:dense1} to a more general
setup.

%---------------------------------------------------------------------------------

%-------------------------------------------------------------------

\end{document}